\documentclass{article}
\usepackage{latexsym,amsmath,amssymb,amscd}
\usepackage[mathscr]{eucal}
\usepackage[dvips]{color}
\begin{document}

\title{On Adiabatic Oscillations of a Stratified Atmosphere on the Flat Earth }
\author{Tetu Makino \footnote{Professor Emeritus at Yamaguchi University, Japan;  E-mail: makino@yamaguchi-u.ac.jp}}
\date{\today}
\maketitle

\newtheorem{Lemma}{Lemma}
\newtheorem{Proposition}{Proposition}
\newtheorem{Theorem}{Theorem}
\newtheorem{Definition}{Definition}
\newtheorem{Remark}{Remark}
\newtheorem{Corollary}{Corollary}
\newtheorem{Notation}{Notation}
\newtheorem{Assumption}{Assumption}

\numberwithin{equation}{section}

\begin{abstract}
We consider the oscillations of the atmosphere around a stratified back ground density and entropy distribution under the gravitation on the flat Earth. The atmosphere is supposed to be an ideal gas and the motion is supposed to be governed by the compressible Euler equations. The density distribution of the back ground equilibrium is supposed to touch the vacuum at the finite height of the stratosphere. Considering the linearized approximation for small perturbations, we show  that time periodic oscillations with a sequence of time periods which accumulate to infinity, say, slow and slow oscillations, so called `g-modes', can appear when the square of the
Brunt-V\"{a}is\"{a}l\"{a} frequency is positive everywhere for the considered back ground equilibrium. 

{\it Key Words and Phrases.} Euler equations,  Atmospheric Oscillations, Vacuum boundary, Brunt-V\"{a}is\"{a}l\"{a} frequency, Gravity modes. Eigenvalue problem  of Strum-Liouville type.


\end{abstract}

\section{Introduction}

We investigate the motion of an atmosphere on the flat earth under the constant gravitational force. We suppose that the atmosphere consists of an ideal gas with the most simple equation of state and the motion is adiabatic motion governed by the compressible Euler equations.

The pioneering mathematical investigation of the oscillations of an atmosphere can be found as the paper `On the Vibrations of an Atmosphere' by Lord Rayleigh, 1890, \cite{Rayleigh}. He wrote

\begin{quote}
In order to introduce greater precision into our ideas respecting the behavior of the Earth's Atmosphere, it seems advisable to solve any problems that may present themselves, even though the search for simplicity may lead us to stray rather far from the actual question.
It is supposed here to consider the case of an atmosphere composed of gas which obeys Boyle's law, viz. such that the pressure is always proportional to the density. And in the first instance we shall neglect the curvature and rotation of the Earth, supposing that the strata of equal density are parallel planes perpendicular to the direction in which gravity acts.
\end{quote}

Our investigation in this article will be done under the same spirit as that of Lord Rayleigh quoted above, except for his starting point that the back ground state is supposed to be that of the isothermal stratified gas, say, $\rho=\rho_0\exp(-\mathsf{g}z/A)$, where $\rho$ is the density and $z$ is the height, $\rho_0, \mathsf{g}, A $ being positive constants. Instead we are interested in a back ground equilibrium like $\rho=C (z_+-z)^{\nu}$, where $z_+$ is the height of the stratosphere, $C, \nu (>1)$ being positive constants, which touches the vacuum $\rho=0$ 
on $z >z_+$ at the height $z=z_+$. Then the vacuum boundary at which the gas touches the vacuum should be treated as a free boundary, and requires a delicate treatise in the mathematically rigorous view point. 

First of all we consider the problem under the linearized approximation for small perturbations around the back ground density and entropy distributions. The aim of the study should be to clarify the spectral property of the differential operator which governs the perturbations around the back ground stratified density distribution in a suitable functional space in view of the linearized approximation, but, especially we are interested in the existence of a sequence of eigenvalues which accumulate to $0$, namely, the existence of time periodic oscillations whose time periods are long and long. This kind of modes are called `g-modes' in the context of Helioseismology, and generally of astroseismology,  by astrophysicists.
This article tries to prove the existence of such a  sequence of modes in the simple situation of the motion under the constant gravitation over the flat Earth. Such a sequence of eigenvalues accumulating to $0$ cannot appear if the background state is isentrooic and the 
 Brunt-V\"{a}is\"{a}l\"{a}
frequency vanishes everywhere. So the existence of g-modes is an effect of buoyancy.\\

Let us describe the problem more precisely. \\

We consider the motions of the atmosphere on the flat earth governed by the Euler equations
\begin{subequations}
\begin{align}
&\frac{\partial{\rho}}{\partial t}+\sum_{k=1}^3\frac{\partial}{\partial x^k}({\rho} v^k)=0, \label{Ea} \\
&{\rho}\Big(\frac{\partial v^j}{\partial t}+\sum_{k=1}^3 v^k\frac{\partial v^j}{\partial x^k}\Big)
+\frac{\partial P}{\partial x^j}+{\rho}\frac{\partial\Phi}{\partial x^j}=0, \quad j=1,2,3,
 \label{Eb} \\
&\rho\Big(\frac{\partial S}{\partial t}+\sum_{k=1}^3 v^k\frac{\partial S}{\partial x^k}\Big)=0, \label{Ec} \\
&\Phi =\mathsf{g} x^3. \label{Ed}
\end{align}
\end{subequations} 
Here $t\geq 0, \mathbf{x}=(x^1,x^2,x^3) \in \Omega:=\{\mathbf{x}\in\mathbb{R}^3 | 
x^3 >0\}$. The unknowns ${\rho} \geq 0, P, S$ are mass density, pressure, entropy density,
and
 $\displaystyle \mathbf{v}=v^1\frac{\partial}{\partial x^1}+v^2\frac{\partial}{\partial x^2}+v^3\frac{\partial}{\partial x^3}$ is the velocity field. $\mathsf{g}$ is a positive constant. The boundary condition is
\begin{equation}
{\rho} v^3=0 \quad\mbox{on}\quad x^3=0, \label{B}
\end{equation}
and the initial condition is 
\begin{equation}
\rho=\overset{\circ}{\rho}(\mathbf{x}),\quad
\mathbf{v}=\overset{\circ}{\mathbf{v}}(\mathbf{x})=
(\overset{\circ}{v}^1(\mathbf{x}),
\overset{\circ}{v}^2(\mathbf{x}),
\overset{\circ}{v}^3(\mathbf{x}))\quad\mbox{at}\quad t=0. \label{I}
\end{equation}

We assume that $P$ is a function of $(\rho, S)$, and put the following

\begin{Assumption}
$P$ as the function of $(\rho, S)$ is  given by
\begin{equation}
P=\rho^{\gamma}\exp\Big(\frac{S}{\mathsf{C}_V}\Big) \quad\mbox{for}\quad \rho \geq 0, \label{ES}
\end{equation}
where $\gamma$ and $\mathsf{C}_V$ are positive constant such that $1<\gamma <2$. 
\end{Assumption}

In this article we denote 
\begin{equation}
\nu:=\frac{1}{\gamma-1}. \label{Nu}
\end{equation}\\

Let us fix a stratified equilibrium $\rho=\bar{\rho}, P=\bar{P}, S=\bar{S}$, which are functions of $x^3$ only such that 
\begin{equation}
\{\bar{\rho} >0 \}=\Pi:=\{ 0 \leq z < z_+\}.
\end{equation}

Here and hereafter we denote $x=x^1, y=x^2, z=x^3$.

We consider the Eulerian perturbations $$\mbox{\boldmath$\xi$}=\sum \xi^k\frac{\delta}{\partial x^k}=
\delta\mathbf{x}=\sum \delta x^k\frac{\partial}{\partial x^k},$$
$ \delta\rho, \delta P, \delta S$ at this fixed equilibrium. Here we use the Lagrangian co-ordinates which will be denoted by the diversion of the letter $\mathbf{x}=(x^1,x^2,x^3)=(x,y,z)$ of the Eulerian co-ordinates. \\

Here let us recall the definition of the Euler perturbation $\delta Q$ and the Lagrange perturbation 
$\Delta Q$ of a quantity $Q$:
\begin{align*}
&\Delta Q(t,\mathbf{x})=Q(t, \mbox{\boldmath$\varphi$}(t,\mathbf{x}))-\bar{Q}(\mathbf{x}), \\
&\delta Q(t,\mathbf{x})=Q(t,\mbox{\boldmath$\varphi$}(t,\mathbf{x}))-\bar{Q}
(\mbox{\boldmath$\varphi$}(t,\mathbf{x})),
\end{align*}
where $\underline{\mathbf{x}}=\mbox{\boldmath$\varphi$}(t,\mathbf{x})=\mathbf{x}+\mbox{\boldmath$\xi$}(t,\mathbf{x})$ is the steam line given by 
$$\frac{\partial}{\partial t}\mbox{\boldmath$\varphi$}(t,\mathbf{x})=\mathbf{v}(t,\mbox{\boldmath$\varphi$}(t,\mathbf{x})),
\quad \mbox{\boldmath$\varphi$}(0,\mathbf{x})=\mathbf{x}. $$ \\

We assume that $\overset{\circ}{\rho}=\bar{\rho}$, or, the initial Lagrangian perturbation
$\Delta\rho|_{t=0}=0$, so that
\begin{equation}
\rho(t,\mathbf{x}+\mbox{\boldmath$\xi$})=\bar{\rho}(\mathbf{x})+\Delta\rho(t,\mathbf{x} )=
\frac{\bar{\rho}(\mathbf{x})}{\mathrm{det}J(t,\mathbf{x})},
\end{equation}
where
\begin{equation}
J(t,\mathbf{x}):= \Big(\delta_k^j+\frac{\partial\xi^j}{\partial x^k}\Big)_{j,k}.
\end{equation}
We suppose $\Delta S|_{t=0}=0$, so that $S(t,\mathbf{x}+\mbox{\boldmath$\xi$}(t, \mathbf{x}))=\bar{S}(\mathbf{x})$.

So $\mathbf{x}$ runs over the fixed domain $\Pi=\{ 0\leq z < z_+\}$. \\

The linearized approximation of the equation which governs the perturbations 
turns out to be
\begin{equation}
\frac{\partial^2\mbox{\boldmath$\xi$}}{\partial t^2}+\mathbf{L}\mbox{\boldmath$\xi$}=0, \label{WE}
\end{equation}
where 
\begin{align}
&\mathbf{L}\mbox{\boldmath$\xi$}=\frac{1}{\bar{\rho}}\mathrm{grad}\delta P+\frac{\mathsf{g}}{\bar{\rho}}\delta\rho\mathbf{e}_3, \\
&\delta\rho=-\mathrm{div}(\bar{\rho}\mbox{\boldmath$\xi$}), \\
&\delta P=\overline{\frac{\gamma P}{\rho}}\delta\rho+\gamma\mathscr{A}\bar{P}\xi^3.
\end{align}
Here $\mathbf{e}_3$ means the unit vector $\partial/\partial x^3=\partial/\partial z$ 
 Hereafter we shall denote
$$\mathbf{e}_1=\frac{\partial}{\partial x^1}=\frac{\partial}{\partial x},\quad
\mathbf{e}_2=\frac{\partial}{\partial x^2}=\frac{\partial}{\partial y},\quad
\mathbf{e}_1=\frac{\partial}{\partial x^3}=\frac{\partial}{\partial z}.
$$

We define 
\begin{equation}
\mathscr{A}:=\overline{\frac{1}{\rho}\frac{d\rho}{dz}}
-\overline{\frac{1}{\gamma P}\frac{dP}{dz}}
=-\frac{1}{\gamma\mathsf{C}_V}\overline{\frac{dS}{dz}}.
\end{equation}

We shall denote
\begin{align}
&c^2:=\overline{\frac{\gamma P}{\rho}}, \\
&\mathscr{N}^2:=-\mathsf{g}\mathscr{A}=-\mathsf{g}
\Big(\overline{\frac{1}{\rho}\frac{d\rho}{dz}}+\frac{\mathsf{g}}{c^2}\Big).
\end{align}

In the physical context the quantity  $c^2$ is the square of the sound speed and $\mathscr{N}^2$ is the square of the  Brunt-V\"{a}is\"{a}l\"{a} frequency.

Then we see 
\begin{equation}
\gamma\mathscr{A}\bar{P}=c^2\overline{\frac{d\rho}{dz}}+\mathsf{g}\bar{\rho}, 
\end{equation}
therefore
\begin{equation}
\delta P=c^2\delta\rho+\Big(c^2\overline{\frac{d\rho}{dz}}+\mathsf{g}\bar{\rho}\Big)\xi^3.
\end{equation}

The boundary condition is
\begin{equation}
\xi^3=(\mbox{\boldmath$\xi$}|\mathbf{e}_3)=0 \quad\mbox{on}\quad z=0. \label{BCxi}
\end{equation}

The initial condition is
\begin{equation}
\mbox{\boldmath$\xi$}=0,\quad
\frac{\partial\mbox{\boldmath$\xi$}}{\partial t}=\overset{\circ}{\mathbf{v}}
\quad \mbox{at}\quad t=0.
\end{equation}

In this article we investigate time periodic solutions of the equation \eqref{WE}. This is the eigenvalue problem associated with the operator $\mathbf{L}$. In the view point of the functional analysis a  clarification of the spectral property of the operator $\mathbf{L}$ is desired, but it is not yet done completely. We shall prove the existence of a sequence of eigenvalues which accumulate to $0$ under 
the assumption that $\inf_{0<z<z_+}\mathscr{N}^2 >0$ (Assumption \ref{PN})
on the value distributions of 
the  square of the Brunt-V\"{a}is\"{a}l\"{a} frequency
$\mathscr{N}^2$.

\section{Equilibrium for a prescribed entropy distribution}

In this section we establish the existence of  equilibria which enjoy good properties used in the following consideration on $\mathbf{L}$.\\

We put the following 

\begin{Definition}
A stratified equilibrium $(\bar{\rho}, \bar{S})$ is said to be admissible if 

1) 
\begin{equation}
\{\bar{\rho} >0 \}=\Pi:=\{ 0 \leq z < z_+\}.
\end{equation}

2)
 $$ z \mapsto \bar{\rho}, \bar{S} \in C^{\infty}(\Pi);$$

3)

$$ \frac{d\rho}{dz} < 0\quad\mbox{for}\quad 0\leq z <z_+, $$

4) 
\begin{equation}
 \bar{\rho} =C_{\rho}(z_+-z)^{\nu}(1+[z_+-z]_1) \quad \mbox{for}\quad 0<z_+-z\ll 1 
\end{equation}
and $z \mapsto \bar{S}$ is analytic at $z_+$, that is,
\begin{equation}
\bar{S}=[z_+-z]_0 \quad\mbox{for}\quad |z_+-z|\ll 1.
\end{equation}
Here $C_{\rho}$ is a positive constant and
we use the following notation.
\end{Definition}

\begin{Notation}
$[X]_K$ stands for various convergent power series of the form
$\displaystyle \sum_{k\geq K} a_kX^k$

\end{Notation}

We claim

\begin{Theorem}\label{Th.1}
 Let a smooth function $\Sigma$ on $\mathbb{R}$ and a positive number 
$z_+$  be given. Assume that it holds, for $\eta>0$, that
\begin{equation}
\gamma +\frac{\gamma-1}{\mathsf{C}_V}\eta\frac{d}{d\eta}\Sigma(\eta) >0. \label{pDP}
\end{equation}
Then there exists an
admissible  equilibrium
$(\bar{\rho}, \bar{S})$ such that $\bar{S}=\Sigma(\bar{\rho}^{\gamma-1})$ .
\end{Theorem}

Proof . Consider the functions $f^P, f^u$ defined by 
\begin{align}
&f^P(\rho):=\rho^{\gamma}\exp\Big[\frac{\Sigma(\rho^{\gamma-1})}{\mathsf{C}_V}\Big], \\
&f^u(\rho):=\int_0^{\rho}
\frac{Df^P(\rho')}{\rho'}d\rho'
\end{align}
for $\rho >0$. Thanks to the assumption \eqref{pDP} we have
$$Df^P(\rho) >0$$
for $\rho >0$, and there exists a smooth function $\Lambda$ on $\mathbb{R}$ such that
$\Lambda(0)=0$ and
\begin{equation}
f^P(\rho)=\mathsf{A}\rho^{\gamma}(1+\Lambda(\rho^{\gamma-1}))
\end{equation}
for $\rho >0$. Here $\mathsf{A}:=\exp(\Sigma(0)/\mathsf{C}_V)$ is a positive constant. 
Then we have
\begin{equation}
u=f^u(\rho)=\frac{\gamma\mathsf{A}}{\gamma-1}\rho^{\gamma-1}(1+\Lambda_u(\rho^{\gamma-1}))
\end{equation}
for $\rho >0$, where $\Lambda_u$ is a smooth function on $\mathbb{R}$ such that 
$\Lambda_u(0)=0$, and  the inverse function $f^{\rho}$ of $f^u$ 
\begin{equation}
f^{\rho}(u)=\Big(\frac{\gamma-1}{\gamma\mathsf{A}}\Big)^{\frac{1}{\gamma-1}}
(u \vee 0)^{\frac{1}{\gamma-1}}(1+\Lambda_{\rho}(u))
\end{equation}
is given so that
$ \rho=f^{\rho}(u) \Leftrightarrow u=f^u(\rho)$ for $u>0 (\rho >0)$. Here $u\vee 0$ stands for $\max( u, 0)$ and  $ \Lambda_{\rho}$ are smooth functions on
$\mathbb{R}$ such that $\Lambda_{\rho}(0)=0$.

Therefore the problem is reduced to that for barotropic case to solve
$$ \frac{du}{dz}+\mathsf{g}=0, $$
which has the solution 
$$u=\mathsf{g}(z_+-z).$$
Of course $u >0 \Leftrightarrow z<z+$ and we are putting 
\begin{align}
 \rho&=f^{\rho}(\mathsf{g}(z_+-z)) \nonumber \\
&=\Big(\frac{\gamma-1}{\gamma\mathsf{A}}\Big)^{\frac{1}{\gamma-1}}
\mathsf{g}((z_+-z) \vee 0)^{\frac{1}{\gamma-1}}(1+\Lambda_{\rho}(\mathsf{g}(z_+-z))
\end{align}
This is the required admissible equilibrium.
$\square$\\

Hereafter in this article we fix an admissible  equilibrium $(\bar{\rho}, \bar{S})$.

\section{Self-adjoint realization of $\mathbf{L}$}

In this article we consider perturbations which are periodic in $x$- and $y$-coordinates. Therefore $\Pi$ will denotes the space $(\mathbb{R}/x_+\mathbb{Z}) \times
(\mathbb{R}/y_+\mathbb{Z}) \times [0, z_+[$, that is, a function $f$ on $\Pi$ is such that
$$ f(x+lx_+, y+my_+,  z) = f(x,y,z)\quad \mbox{for}\quad \forall l,m \in \mathbb{Z}.$$
Here the periods $x_+, y_+$ are arbitrarily fixed to be positive numbers.
\\

We are considering the differential operator
\begin{equation}
\mathbf{L}\mbox{\boldmath$\xi$}=\frac{1}{\rho}\mathrm{grad}\delta P
+\frac{\mathsf{g}}{\rho}\delta\rho \mathbf{e}_3
\end{equation}
where
\begin{subequations}
\begin{align}
\delta\rho&=-\mathrm{div}(\rho\mbox{\boldmath$\xi$}), \\ 
\delta P&=\frac{\gamma P}{\rho}\delta \rho+\gamma \mathscr{A}P\xi^3.
\end{align}
\end{subequations}

 Here and hereafter the bars to denote the quantities evaluated at the fixed equilibrium are omitted, that is, $\rho, P$ etc stand for $\bar{\rho}, \bar{P}$ etc\\

Let us consider the operator $\mathbf{L}$ in the Hilbert space $\mathfrak{H}=L^2(\Pi, \rho dx), \mathbb{C}^3)$ endowed with the norm $\|\mbox{\boldmath$\xi$}\|_{\mathfrak{H}}$ defined by
\begin{equation}
\|\mbox{\boldmath$\xi$}\|_{\mathfrak{H}}^2=\int_{\Pi}
|\mbox{\boldmath$\xi$}(\mathbf{x})|^2
\rho(\mathbf{x})d\mathbf{x}
=\int_0^{x_+}\int_0^{y_+}\int_0^{z_+}|\mbox{\boldmath$\xi$}(\mathbf{x})|^2
\rho(z)
dzdydx.
\end{equation}

We shall use

\begin{Notation}
For complex number $z=x+\sqrt{-1}y, x, y \in \mathbb{R}$, the complex conjugate is denoted by
$z^*=x-\sqrt{-1}y$. Thus, for $\displaystyle \mbox{\boldmath$\xi$}=\sum \xi^k\frac{\partial}{\partial x^k}, \xi^k \in
\mathbb{C}$, we denote $\displaystyle \mbox{\boldmath$\xi$}^*=\sum (\xi^k)^*\frac{\partial}{\partial x^k}$.
\end{Notation}

First we observe $\mathbf{L}$ restricted on $C_0^{\infty}(\Pi^o, \mathbb{C}^3)$. Here and hereafter we denote
\begin{equation}
\Pi^o:=(\mathbb{R}/x_+\mathbb{Z})\times(\mathbb{R}/y_+\mathbb{Z})\times ]0,z_+[
=\Pi\setminus \{z=0\}.
\end{equation}

We look at
\begin{align}
\mathbf{L}\mbox{\boldmath$\xi$}&=
\mathrm{grad}\Big(-\frac{\gamma P}{\rho^2}\mathrm{div}(\rho\mbox{\boldmath$\xi$})+
\frac{\gamma\mathscr{A}P}{\rho}\xi^3\Big) + \nonumber \\
&+\frac{\gamma\mathscr{A}P}{\rho^2}\Big(-\mathrm{div}(\rho\mbox{\boldmath$\xi$})+
\frac{d\rho}{dz}\xi^3\Big)\mathbf{e}_3.
\end{align}

Using this expression for $\mbox{\boldmath$\xi$}_{(\mu)} \in C_0^{\infty}(\Pi^o), \mu=1,2$, we have the following formula by integration by parts:
\begin{align*}
(\mathbf{L}\mbox{\boldmath$\xi$}_{(1)}|\mbox{\boldmath$\xi$}_{(2)})_{\mathfrak{H}}&=
\int\frac{\gamma P}{\rho^2}\mathrm{div}(\rho\mbox{\boldmath$\xi$}_{(1)})\mathrm{div}(\rho\mbox{\boldmath$\xi$}_{(2)}^*) + \\
&+\int\frac{\gamma\mathscr{A}P}{\rho}
\Big[(\mbox{\boldmath$\xi$}_{(1)}|\mathbf{e}_3)\cdot\mathrm{div}(\rho\mbox{\boldmath$\xi$}_{(2)}^*)
-\mathrm{div}(\rho\mbox{\boldmath$\xi$}_{(1)})\cdot (\mbox{\boldmath$\xi$}_{(2)}|\mathbf{e}_3)^*\Big] + \\
&+\int\frac{\gamma\mathscr{A}P}{\rho}\frac{d\rho}{dz}(\mbox{\boldmath$\xi$}_{(1)}|\mathbf{e}_3)
(\mbox{\boldmath$\xi$}_{(2)}|\mathbf{e}_3)^*.
\end{align*}

We see that
$\mathbf{L}$ restricted on $C_0^{\infty}(\Pi^o)$ is a symmetric operator.
Of course $C_0^{\infty}(\Pi^o)$ is dense in $\mathfrak{H}$.

Moreover we have
\begin{align*}
(\mathbf{L}\mbox{\boldmath$\xi$}|\mbox{\boldmath$\xi$})&=
\int\frac{\gamma P}{\rho^2}|\mathrm{div}(\rho\mbox{\boldmath$\xi$})|^2 + \\
&+2\sqrt{-1}\mathfrak{Im}\Big[
\int\frac{\gamma\mathscr{A}P}{\rho}\xi^3\cdot\mathrm{div}(\rho\mbox{\boldmath$\xi$}^*)\Big] +
\int\frac{\gamma\mathscr{A}P}{\rho}\frac{d\rho}{dz}|\xi^3|^2.
\end{align*}

Since $\mathscr{A} \in C^{\infty}(\overline{\Pi})$, we have
$$|\mathscr{A}|\sqrt{\frac{\gamma P}{\rho}}\leq C_1$$
on $0<z < z_+$, for $P/\rho =O(z_+-z)$. Therefore
\begin{align*}
\Big|\int
\frac{\gamma\mathscr{A}P}{\rho}\xi^3
\mathrm{div}(\rho\mbox{\boldmath$\xi$})^*\Big|
&\leq C_1\int
\sqrt{ \frac{\gamma P}{\rho} }|\xi^3|
|\mathrm{div}(\rho\mbox{\boldmath$\xi$})| \\
&\leq \frac{C_1}{2}\Big[
\frac{1}{\epsilon}\int \rho |\xi^3|^2+
\epsilon\int
\frac{\gamma P}{\rho^2}|\mathrm{div}(\rho\mbox{\boldmath$\xi$})|^2\Big] \\
&\leq \frac{C_1}{2}\Big[
\frac{1}{\epsilon}\|\mbox{\boldmath$\xi$}\|_{\mathfrak{H}}^2+
\epsilon\int\frac{\gamma P}{\rho^2}|\mathrm{div}\rho\mbox{\boldmath$\xi$})|^2\Big].
\end{align*}

Since 
$\displaystyle \frac{P}{\rho}\frac{d\rho }{dz}=O(\rho)$, we have 
$$\Big|\frac{\gamma\mathscr{A}P}{\rho}\frac{d\rho}{dz}\Big|\leq C \rho $$
Therefore we have
$$
\Big|\int\frac{\gamma\mathscr{A}P}{\rho}\frac{d\rho}{dz}|\xi^3|^2\Big|\leq C_2\|\mbox{\boldmath$\xi$}\|_{\mathfrak{H}}^2.
$$

Thus
$$
(\mathbf{L}\mbox{\boldmath$\xi$}|\mbox{\boldmath$\xi$})_{\mathfrak{H}}\geq
\Big(1-\frac{\epsilon C_1}{2}\Big)\int
\frac{\gamma P}{\rho^2}|\mathrm{div}(\rho\mbox{\boldmath$\xi$})|^2 
-\Big(\frac{C_1}{2\epsilon}+C_2\Big)\|\mbox{\boldmath$\xi$}\|_{\mathfrak{H}}^2. $$
Taking $\epsilon$ so small that $1-\frac{\epsilon C_1}{2} \geq 0$, we get
$$
(\mathbf{L}\mbox{\boldmath$\xi$}|\mbox{\boldmath$\xi$})_{\mathfrak{H}}\geq
-\Big(\frac{C_1}{2\epsilon}+C_2\Big)\|\mbox{\boldmath$\xi$}\|_{\mathfrak{H}}^2. $$

Summing up, $\mathbf{L}$ is bounded from below in $\mathfrak{H}$. Therefore, thanks to \cite[Chapter VI, Section 2.3]{Kato}, we have

\begin{Theorem}
The differential operator $\mathbf{L}$ on $C_0^{\infty}(\Pi^o, \mathbb{C}^3)$ admits the Friedrichs extension, which is a self-adjoint operator, in $\mathfrak{H}$.
\end{Theorem}

 Hereafter we shall denote by the same letter $\mathbf{L}$ this Friedrichs extension. The domain $\mathsf{D}(\mathbf{L})$ of the self-adjoint operator $\mathbf{L}$ is characterized as 
$$\mathsf{D}(\mathbf{L})=\{ \mbox{\boldmath$\xi$} \in \mathfrak{H}_1 | \mathbf{L}\mbox{\boldmath$\xi$} \in \mathfrak{H} \quad\mbox{in distribution sense} \}.$$
Here $\mathfrak{H}_1$  is the closure of $C_0^{\infty}(\Pi^o;\mathbb{C})$ in the Hilbert space $\mathfrak{G}$ endowed with that norm
\begin{align*}
\|\mbox{\boldmath$\xi$}\|_{\mathfrak{G}}^2&=\int\frac{\gamma P}{\rho}|\mathrm{div}(\rho\mbox{\boldmath$\xi$}|^2+
2\sqrt{-1}\mathfrak{Im} \Big[\int \frac{\gamma \mathscr{A}P}{\rho}\xi^3\cdot\mathrm{div}(\rho\mbox{\boldmath$\xi$}^*)\Big]+ \\
&+\int\frac{\gamma \mathscr{A}P}{\rho}\frac{d\rho}{dz}|\xi^3|^2+
K_0\|\mbox{\boldmath$\xi$}\|_{\mathfrak{H}}^2,
\end{align*}
where $K_0$ is a sufficiently large positive number.

We want to clarify the spectral property of the self-adjoint operator $\mathbf{L}$.\\

When $\bar{S}$ is constant, say, the equilibrium is isentropic, so that $\mathscr{A}=0$
and $\mathscr{N}^2=0$  everywhere, the spectral property of $\mathbf{L}$ is clear in some sense. Actually in this case the operator $\mathbf{L}$ turns out to be
$$\mathbf{L}\mbox{\boldmath$\xi$} =-\mathrm{grad}\Big(\frac{c^2}{\rho}g\Big) \qquad\mbox{with}\quad g:=\mathrm{div}(\rho\mbox{\boldmath$\xi$}).
$$
Therefore the vector wave equation
$$\frac{\partial^2\mbox{\boldmath$\xi$}}{\partial t^2}+\mathbf{L}\mbox{\boldmath$\xi$} =0 $$
can be reduced to the scalar wave equation
$$\frac{\partial^2g}{\partial t^2}-\mathrm{div}\Big(\rho\mathrm{grad}\Big(\frac{c^2}{\rho}g\Big)\Big)=0.
$$
Then by the argument developed in \cite{JJTM} we can claim the following :\\

{\it 
Suppose that $\bar{S}$ is constant so that $\mathscr{N}^2=0$ everywhere. Then $\mathbf{L}$ can be considered as a self-adjoint operator on the Hilbert space
\begin{equation}
\mathfrak{F}:=\mathfrak{H}\cap\Big\{g=\mathrm{div}(\rho\mbox{\boldmath$\xi$}) \in
L^2(\Pi, \frac{c^2}{\rho} d\mathbf{x}) \Big| \int_{\Pi}g=0\Big\},
\end{equation}
which is dense in $\mathfrak{H}$, and the spectrum of the operator $\mathbf{L}$ in $\mathfrak{F}$ consists of the eigenvalue $0$ of infinite multiplicity and
sequence of eigenvalues $(\lambda_n)_{n\in\mathbb{N}}$ of finite multiplicity such that
$\lambda_n \rightarrow +\infty$ as $n \rightarrow \infty$.}\\

But we are interested in the case in which $\mathscr{N}^2$ does not vanish identically.
We shall show that if $\mathscr{N}^2 > 0$ everywhere, there can appear sequence of eigenvalues which accumulates to $0$. \\

\section{A class of particular solutions}

Let us consider  perturbations $\mbox{\boldmath$\xi$}$ of  solutions of the particular form
\begin{equation}
\mbox{\boldmath$\xi$}=u(z)\sin lx \mathbf{e}_1+w(z)\cos lx \mathbf{e}_3. \label{401}
\end{equation}
Here the frequency $l$ in $x$ is such that $lx_+ \in 2\pi\mathbb{Z}$. Then the eigenvalue problem 
\begin{equation}
\mathbf{L}\mbox{\boldmath$\xi$}=\lambda\mbox{\boldmath$\xi$}
\end{equation}
turns out to be 
\begin{equation} 
\vec{L}_l\vec{w}=\lambda\vec{w} \label{4002}
\end{equation}
for 
\begin{equation}
\vec{w}=
\begin{bmatrix}
u \\
\\
w
\end{bmatrix},
\end{equation}
where
\begin{equation}
\vec{L}_l\vec{w}=
\begin{bmatrix}
L_l^u \\
\\
L_l^w
\end{bmatrix}
=
\begin{bmatrix}
\displaystyle -\frac{l}{\rho}(\delta\check{P})_l, \\
\\
\displaystyle \frac{1}{\rho}\frac{d}{dz}(\delta\check{P})_l+\frac{\mathsf{g}}{\rho}(\delta\check{\rho})_l,
\end{bmatrix}
\end{equation}
and
\begin{align}
(\delta\check{\rho})_l&=-l\rho u -\frac{d}{dz}(\rho w), \\
(\delta\check{P})_l&=c^2(\delta\check{\rho})_l+
\Big(c^2\frac
{d\rho}{dz}+\mathsf{g}\rho\Big)w \nonumber \\
&=-c^2\rho\Big(lu +\frac{dw}{dz}\Big)+\mathsf{g}\rho w
\end{align}\\

Note that {\it if frequencies $l^{(1)}, l^{(2)}$ such that $l^{(1)}x_+ \in 2\pi\mathbb{Z}, l^{(2)}y_+\in 2\pi\mathbb{Z}$ are given, and if $\displaystyle \vec{w}^{(j)}=
\begin{bmatrix}
u^{(j)} \\
w^{(j)}
\end{bmatrix}
, j=1,2,$ satisfy
$$\vec{L}_{l^{(j)}}\vec{w}^{(j)}=\lambda^{(j)}\vec{w}^{(j)},\qquad j=1,2,
$$
then the perturbation
\begin{align*}
\mbox{\boldmath$\xi$}&=u^{(1)}(z)\sin (\sqrt{\lambda^{(1)}}t)\sin(l^{(1)}x)\mathbf{e}_1+
u^{(2)}(z)\sin (\sqrt{\lambda^{(2)}}t)\sin(l^{(2)}y)\mathbf{e}_2+ \\
&+
[w^{(1)}(z)\sin (\sqrt{\lambda^{(1)}}t)\cos(l^{(1)}x)+
w^{(2)}(z)\sin (\sqrt{\lambda^{(2)}}t)\cos(l^{(2)}y)]\mathbf{e}_3
\end{align*}
solves the wave equation 
$$
\frac{\partial^2\mbox{\boldmath$\xi$}}{\partial t^2}+\mathbf{L}\mbox{\boldmath$\xi$}=0.
$$
Of course, any linear combination with constant coefficients of finite number of perturbations of the above form with various $l^{(1)}, l^{(2)}, \lambda^{(1)}, \lambda^{(2)}$ satisfies the wave equation.}
When $l^{(1)}=0$ or $l^{(2)}=0$, we neglect $u^{(1)}$ or $u^{(2)}$ respectively, since $L_0^u=0$.\\

Keeping in mind this situation, we are going to consider the particular perturbation \eqref{401} and the eigenvalue problem \eqref{4002} without loss of generality.

\subsection{ Case $l=0$ }

Let us consider the case $l=0$, namely, let us consider merely vertical perturbations.

By neglecting $u$ and eliminating $(\delta\check{\rho})_0$, the problem is reduced to
the system
\begin{subequations}
\begin{align}
&\frac{dw}{dz}-\frac{\mathsf{g}}{c^2}w+
\frac{(\delta\check{P})_0}{c^2{\rho}}=0, \label{0l0a} \\
&\frac{d}{dz}(\delta\check{P})_0+\frac{\mathsf{g}}{c^2}(\delta\check{P})_0+
(\mathscr{N}^2-\lambda){\rho}w=0, \label{0l0b}
\end{align}
\end{subequations}
and it is reduced to the single equation
\begin{equation}
-\frac{d}{dz}\Big(c^2\rho\frac{dw}{dz}\Big)
=\lambda\rho w. \label{0l1}
\end{equation}
Let us perform the Liouville transformation. See, e.g., \cite[p.275, Theorem 6]{BirkhoffR} or \cite[p.110]{Yosida}. We put
\begin{align}
&\zeta=\int_0^z\sqrt{\frac{a}{\kappa}}dz, \qquad W=(a\kappa)^{\frac{1}{4}}w, \nonumber  \\
&q=
\frac{1}{4}\frac{a}{\kappa}\Big[\frac{d^2}{dz^2}\log(a\kappa)-
\frac{1}{4}\Big(\log(a\kappa)\Big)^2+ \nonumber \\
&+\Big(\frac{d}{dz}\log a\Big)\Big(\frac{d}{dz}\log (a\kappa)\Big)\Big]
\end{align}
for
\begin{equation}
a=c^2\rho, 
\qquad \kappa=\rho,
\end{equation}
which transform \eqref{0l1}
to
\begin{equation}
-\frac{d^2W}{d\zeta^2}+qW=\lambda W.
\end{equation}
We see
$$\zeta=\int_0^z\frac{1}{c}dz $$
so that
\begin{align}
&\zeta_+:=\int_0^{z_+}\frac{1}{c}dz < \infty, \\
&\zeta_+-\zeta=2\sqrt{\frac{\nu}{\mathsf{g}}}\sqrt{z_+-z}
(1+[z_+-z]_1)
\quad\mbox{for}\quad 0<z_+-z \ll1,
\end{align}
since
$$c^2=\frac{\mathsf{g}}{\nu}(z_+-z)(1+[z_+-z]_1).
$$
By a tedious calculation we see
\begin{equation}
q=\frac{(2\nu+1)(2\nu-1)}{4}\frac{1}{(\zeta_+-\zeta)^2}(1+[(\zeta_+-\zeta)^2]_1).
\end{equation}
Note that $(2\nu+1)(2\nu-1)/4 >3/4$ for $\nu >1$. Therefore we can claim

{\it 
The operator $\displaystyle -\frac{d^2}{d\zeta^2}+q$ defined on $C_0^{\infty}(]0,\zeta_+[)$ admits he Friedrichs extension, which is a self-adjoint operator, in the Hilbert space $L^2([0,\zeta_+])$. Its spectrum consists of simple eigenvalues
$$\lambda_1<\lambda_2< \cdots < \lambda_n<\cdots \rightarrow +\infty.
$$
}
 Here we note that $\lambda_1>0$. In fact, if $w_1$ is an eigenfunction, then 
$$\int c^2\rho\Big|\frac{dw_1}{dz}\Big|^2dz =\lambda_1\int \rho |w_1|^2dz $$
implies $\lambda_1\geq0$; If $\lambda_1=0$, then $w=\mbox{Const.}$ so that $w=0$ thanks to the boundary condition $w(0)=0$, a contradiction; Hence $\lambda_1>0$.

Clearly these eigenvalues are those of $L_0^w$ in $L^2([0,z_+], \rho(z)dz)$  and of $\mathbf{L}$
in $\mathfrak{H}$. Therefore we can claim

\begin{Theorem}\label{Efl0}
There exists a sequence of eigenvalues $(\lambda^{(0)}_n)_{n\in\mathbb{N}}$ of $\mathbf{L}$ such that
$$ \lambda_1^{(0)} < \lambda_2^{(0)}<\cdots <  \lambda_n^{(0)}<\cdots \rightarrow +\infty.$$ They are simple eigenvalues of $L_0^w$, or of $\vec{L}_0$ restricted onto $\{ u=0\}$.
\end{Theorem}

\subsection{ Case $l >0$ }

Let us  consider the equations when $l \not=0$. We can suppose $ l >0$ without loss of generality, since, otherwise, we can consider $(-u, w)$ instead of $(u, w)$, which transformation reduces $\vec{L}_l$ to $\vec{L}_{-l}$. \\

Corresponding to $\mathfrak{H}$, we are considering the Hilbert space $\mathfrak{X}_l=L^2([0,z_+], \rho(z)dz); \mathbb{C}^2)$ of functions $\vec{w}=(u,w)^{\top}$ endowed with the norm
$\|\cdot\|_{\mathfrak{X}_l}$ given by
$$\|\vec{w}\|_{\mathfrak{X}_l}^2=\int_0^{z_+}(l^2|u|^2+|w|^2)\rho(z)dz.
$$

The operators $ \vec{L}_l=(L_l^u, L_l^w)^{\top}$ is considered as a self-adjoint operator 
in $\mathfrak{X}_l$ and the domain $\mathsf{D}(\vec{L}_l)$ is
$$\mathsf{D}(\vec{L}_l)=
\overset{\circ}{\mathfrak{W}}_l\cap
\{ \vec{w} | \vec{L}_l\vec{w} \in \mathfrak{X}_l\}.
$$
Here $\overset{\circ}{\mathfrak{W}}_l$ is the closure of $C_0^{\infty}(]0, z_+[)$
in the Hilbert space $\mathfrak{W}_l$ endowed with norm
$\|\cdot\|_{\mathfrak{W}_l}$ defined by 
$$\|\vec{w}\|_{\mathfrak{W_l}}^2=
\|\vec{w}\|_{\mathfrak{X}_l}^2+\|(\delta\check{\rho})_l\|_{L^2(\frac{c^2}{\rho}dz)}^2.$$\\

Of course this domain corresponds to the domain of the Friedrichs extension $\mathbf{L}$,
say, $\vec{w} \in \mathsf{D}(\vec{L}_l)$ if and only if $\mbox{\boldmath$\xi$}$ specified by \eqref{401}
belongs to $\mathsf{D}(\mathbf{L})$.\\

 Consider $\mbox{\boldmath$\xi$}$ of the form \eqref{401}, when
$$\mathrm{div}(\rho\mbox{\boldmath$\xi$})=(l\rho u(z)+\frac{d}{dz}(\rho w))\cos lx=-(\delta\check{\rho})_l,
$$
Then $\|(\delta\check{\rho})_l\|_{L^2(c^2dz/\rho)} < \infty$ means
$$\int\frac{c^2}{\rho}\Big|l\rho u+\frac{d}{dz}(\rho w)\Big|^2dz <\infty.$$
On the other hand
$$\int\frac{c^2}{\rho}|l\rho u|^2dz =l^2\int c^2\rho|u|^2dz <\infty$$
for $\vec{w} \in \mathfrak{X}_l$. Therefore
$$\int\frac{c^2}{\rho}\Big|\frac{d}{dz}(\rho w)\Big|^2dz <\infty.$$
Since $\inf_{0<z<z_+/2}c^2/\rho >0$, we see 
$$C=\Big[\int_0^{z_+/2}\Big|\frac{d}{dz}(\rho w)\Big|^2dz \Big]^{1/2} <\infty.$$
Since $\vec{w}$ belongs to the closure of $C_0^{\infty}(]0,z_+[)$, we can claim
$$\rho(z)w(z) \leq C\sqrt{z}\qquad\mbox{for}\quad 0<z\leq z_+/2.$$
Since $\inf_{0<z<z_+/2}\rho(z) >0$, we have
$w(z)=O(\sqrt{z})$ as $z\rightarrow 0$, and the boundary condition \eqref{BCxi} at $z=0$ is satisfied.

We can claim

\begin{Proposition} \label{Prop.BC}
Let $\vec{w}=(u, w)^{\top} \in \mathfrak{W}_l$. Then 
$\vec{w} \in \overset{\circ}{\mathfrak{W}}_l $ if
\begin{align} 
& w=0 \quad\mbox{at}\quad z=0, \nonumber \\
\mbox{and} & \quad 
w=O((z_+-z)^{\alpha})\quad\mbox{with}\quad \alpha >-\frac{\nu}{2}\quad\mbox{as}\quad z \rightarrow z_+-0.\label{wBC}
\end{align}
\end{Proposition}

Proof can be done by taking a sequence of cut-off functions $(\varphi_n)_{n=1,2,\cdots}$ in $C_0^{\infty}(\mathbb{R})$ such that $\varphi_n(z)=1$ for $ z\leq z_+-\frac{1}{n}$, $\varphi_n(z)=0$ for $z_+-\frac{1}{2n} \leq z$, $0 \leq \varphi_n(z) \leq 1, |d\varphi_n/dz|\leq Cn$ and considering $\varphi_n\cdot \vec{w}$ for $\vec{w} \in \mathfrak{W}_l$. Let us omit the details. \\

First we note the following

\begin{Proposition}
Let $l >0$. 1) If $\mathscr{A}=0$ identically, then $\lambda=0$ is an eigenvalue of
$\vec{L}_l$ with infinite multiplicity. 2) If $\mathscr{A}\not= 0$ almost everywhere, then
$0$ is not an eigenvalue of $\vec{L}_l$, that is, $\mathrm{Ker}\vec{L}_l =\{ \vec{0}\} $.
\end{Proposition}

Proof. 1) Suppose $\mathscr{A}=0$ identically. Then $(\delta\check{P})_l=c^2
(\delta\check{\rho})_l$ holds and, if they vanish, $\vec{L}_l\vec{w}=0$. It is the case when
\begin{align*}
u&=-\frac{1}{l\rho}\frac{d}{dz}(\rho\Upsilon), \\
w&=\Upsilon(z),
\end{align*}
provided that $\Upsilon \in C_0^{\infty}(]0, z_+[)$. Therefore $\mathrm{dim.Ker}\vec{L}_l=\infty$. 2) Suppose $\mathscr{A}\not=0$ a.e. If $\vec{w} \in \mathrm{Ker}\vec{L}_l$, we have
$(\delta\check{P})_l=c^2(\delta\check{\rho})_l=0$, so that
$$\gamma \mathscr{A}P w=(\delta\check{P})_l-c^2(\delta\check{\rho})_l=0, $$
which implies $w=0$, since $\gamma\mathscr{A}P\not=0$ a.e., and
$$u=-\frac{1}{l\rho}\Big((\delta\check{\rho})_l+\frac{d}{dz}(\rho w)\Big)=0,
$$
that is, $\vec{w}=(0,0)^{\top}$. $\square$ \\

Moreover we try to find eigenvalue $\lambda \not=0$. We want to give a mathematically rigorous justification of the existence of `g-modes', namely, eigenvalues accumulating to $0$ and  `p-modes', namely, eigenvalues accumulating to $+\infty$\\

By eliminating $u$ and $(\delta\check{\rho})_l$, the problem is reduced to
\begin{subequations}
\begin{align}
&\frac{dw}{dz}-\frac{\mathsf{g}}{c^2}w+\Big(1-\frac{l^2c^2}{\lambda}\Big)
\frac{(\delta\check{P})_l}{c^2{\rho}}=0, \label{5a} \\
&\frac{d}{dz}(\delta\check{P})_l+\frac{\mathsf{g}}{c^2}(\delta\check{P})_l+
(\mathscr{N}^2-\lambda){\rho}w=0. \label{5b}
\end{align}
\end{subequations}\\

It would be natural to eliminate $(\delta\check{P})_l$ to deduce a single equation for $w$ as for the case $l=0$. The result would be
\begin{equation}
-\frac{d}{dz}\Big(\frac{c^2\rho}{\mathcal{Q}}\frac{dw}{dz}\Big)+
\Big[\mathsf{g}\frac{\mathcal{Q}}{\rho}\frac{d}{dz}\frac{\rho}{\mathcal{Q}}+
\frac{\mathsf{g}^2}{c^2}(\mathscr{N}^2-\lambda)\Big]\frac{\rho}{\mathcal{Q}}w=0, \label{Elw}
\end{equation}
where
\begin{equation}
\mathcal{Q}:=1-\frac{l^2c^2}{\lambda}.
\end{equation}
But, since $l^2\not=0$ and $c^2 >0$, 
$$ c^2 \sim \frac{\mathsf{g}}{\nu}(z_+-z) \quad\mbox{as}\quad z \rightarrow z_+-0,
$$
for each $\lambda \in ]0, \lambda_0]$, $\lambda_0$ being a small positive number, there exists a unique $z=z(\lambda) \in ]0, z_+[$, namely, the `turning point', such that 
\begin{align*}
& \mathcal{Q}(z) <0 \qquad \mbox{for}\quad 0<z<z(\lambda) \\
&\mathcal{Q}(z(\lambda))=0 \\
&\mathcal{Q}(z) >0 \qquad \mbox{for}\quad z(\lambda)<z<z_+.
\end{align*}
Therefore the coefficients of the equation \eqref{Elw} have poles at the turning point. This may cause trouble during the analysis. So, we should seek another reduction to a single second order equation
which has no singularity in the interior of the interval $[0, z_+]$ and 
 we can treat without trouble. 

\begin{Remark}
We note that the equation \eqref{Elw} is equivalent to \cite[p. 23, (22)]{Yih}, which is nothing but
$$-\lambda\mathcal{Q}c^2\frac{d}{dz}\Big(\rho\frac{dw}{dz}\Big)
-\lambda\mathcal{Q}^2\frac{d}{dz}\Big(\frac{c^2\rho}{\mathcal{Q}})\frac{dw}{dz}
+\lambda\Big[\cdots \Big]\mathcal{Q}\rho c^2w=0,
$$
given by multiplying \eqref{Elw} by $\lambda\mathcal{Q}^2$,
while we can verify
$$\mathcal{Q}^2\frac{d}{dz}\frac{c^2\rho}{\mathcal{Q}}=\frac{dc^2}{dz}\rho.
$$
Actually the symbols $\alpha^2, c_s^2, \sigma^2, (\cdot)'$ in \cite{Yih} read $l^2, c^2, \lambda, \frac{d}{dz}(\cdot)$
respectively in this article. This equation could be said to be a natural generalization of the so called `Taylor-Goldstein equation' which is often used in the study of stratified fluid motions successfully, but it seems not suitable for mathematical investigations here. 
\end{Remark}

In order to avoid the above mentioned trouble, let us introduce the variable $\eta$ instead of $(\delta\check{P})_l$, suggested by the argument by D. O. Gough in \cite{Gough},  by
\begin{align}
\eta&:=(\delta\check{P})_l+\frac{d\bar{P}}{dz}w=(\delta\check{P})_l-\mathsf{g}{\rho}w, \nonumber \\
&=-c^2\rho\Big(lu+\frac{dw}{dz}\Big) \label{Equeta}
\end{align}
which is 
the Lagrangian perturbation $\Delta\check{P}$ in the linearized approximation. Let us note the following

\begin{Proposition} \label{Prop.ueta}
We consider the correspondence between the sets of  variables $(u, w)$ and  $(w, \eta)$ defined by \eqref{Equeta}. Then  $(u,w)^{\top} \in \mathfrak{W}_l$ if and only if
$w \in L^2(\rho dz), \eta \in L^2({dz}/{c^2\rho}) $ and 
\begin{equation}
\frac{dw}{dz}+\frac{\eta}{c^2\rho} \in L^2(\rho dz).
\end{equation}
\end{Proposition}

Then
the set of equations \eqref{5a},\eqref{5b} is equivalent to
\begin{subequations}
\begin{align}
&\frac{dw}{dz}+A_{11}w+A_{12}\eta=0, \label{8a} \\
&\frac{d\eta}{dz}+A_{21}w+A_{22}\eta=0, \label{8b}
\end{align}
\end{subequations}
where 
\begin{align}
& A_{11}=-\frac{l^2\mathsf{g}}{\lambda}, \qquad A_{12}=\Big(1-\frac{l^2c^2}{\lambda}\Big)
\frac{1}{c^2{\rho}}, \nonumber  \\
& A_{21}=\frac{l^2\mathsf{g}^2}{\lambda}\Big(1-\frac{\lambda^2}{l^2\mathsf{g}^2}\Big){\rho}, \qquad A_{22}=\frac{l^2\mathsf{g}}{\lambda}.
\end{align}\\

We see that the set of equations \eqref{8a}, \eqref{8b} reduces to the single equation
\begin{equation}
\frac{d^2\eta}{dz^2}+\frac{1}{\mathsf{H}[{\rho}]}\frac{d\eta}{dz}+
\Big(\frac{\lambda}{c^2}-l^2\Big(1-
\frac{\mathscr{N}^2}{\lambda}\Big)\Big)\eta=0. \label{11}
\end{equation}
Here and hereafter we use the notation
\begin{equation}
\frac{1}{\mathsf{H}[Q]}:=-\frac{d}{dz}\log Q
\end{equation}
for any quantity $Q$ which is a function of $z$. $\mathsf{H}[Q]$ is so called the scale height of $Q$.

In fact, \eqref{11} can be derived by eliminating $w$ from \eqref{8a}\eqref{8b}, provided that $\lambda \not=l\mathsf{g}$, when $A_{21}\not=0$, and can be verified directly when $\lambda=l\mathsf{g}$. 

Note that
\begin{equation}
\mathscr{N}^2=-\frac{\mathsf{g}^2}{c^2}+\frac{\mathsf{g}}{\mathsf{H}[{\rho}]}.
\end{equation}\\

%%%%%%%%%%%%%%%%

Later we shall use the following

\begin{Proposition}\label{Prop.FS}
The system \eqref{8a}\eqref{8b} admits 
a fundamental matrix of solutions $\mbox{\boldmath$\Phi$}_S=[\mbox{\boldmath$\varphi$}_{S1} \quad \mbox{\boldmath$\varphi$}_{S2}]$ of the form
\begin{subequations}
\begin{align}
&\mbox{\boldmath$\varphi$}_{S1}=
\begin{cases}
\begin{bmatrix}
1+[z_+-z]_1 \\
\\
(z_+-z)^{\nu+1}[z_+-z]_0
\end{bmatrix}
\quad\mbox{if}\quad \nu \in \mathbb{N}
\\
\\
\begin{bmatrix}
1+[z_+-z]_1 \\
\\
0
\end{bmatrix}
\quad\mbox{if}\quad \nu \not\in \mathbb{N}, 
\end{cases} \label{FS1a} \\
&\mbox{\boldmath$\varphi$}_{S2}=
\begin{bmatrix}
\displaystyle -\frac{1}{\mathsf{g}C_{\rho}}(z_+-z)^{-\nu}(1+[z_+-z]_1)+[z_+-z]_1 \\
\\
1+[z_+-z]_1
\end{bmatrix}
\label{FS1b}
\end{align}
\end{subequations}
\end{Proposition}

Proof. 
Take the variables
\begin{equation}
s=z_+-z,\qquad \hat{\eta}=s^{-\nu}\eta.
\end{equation}
Then the system \eqref{8a}\eqref{8b} turns out to be
\begin{subequations}
\begin{align}
&s\frac{dw}{ds}=\hat{A}_{11}w+\hat{A}_{12}\hat{\eta},  \label{FS3a}\\
&s\frac{d\hat{\eta}}{ds}=\hat{A}_{21}w+\hat{A}_{22}\hat{\eta},
\end{align}
\end{subequations}
where
\begin{align}
&\hat{A}_{11}=s A_{11}=-\frac{l^2\mathsf{g}}{\lambda}s, \nonumber \\
&\hat{A}_{12}=s^{\nu+1} A_{12}=\frac{\nu}{\mathsf{g}C_{\rho}}(1+[s]_1), \nonumber \\
&\hat{A}_{21}=s^{1-\nu} A_{21}=\frac{l^2\mathsf{g}^2}{\lambda}\Big(
1-\frac{\lambda^2}{l^2\mathsf{g}^2}\Big)C_{\rho}s(1+[s]_1), \nonumber \\
&\hat{A}_{22}=-\nu + s A_{22}=-\nu+\frac{l^2\mathsf{g}}{\lambda}s.
\end{align}
Thus
\begin{equation}
s\frac{d}{ds}
\begin{bmatrix}
w \\
\hat{\eta}
\end{bmatrix}
=(K_0+\sum_{m\geq 1}s^mK_m)
\begin{bmatrix}
w \\
\hat{\eta}
\end{bmatrix}, \label{FS5}
\end{equation} 
where
$$K_0=
\begin{bmatrix}
0 & \frac{\nu}{\mathsf{g}C_{\rho}} \\
0 & -\nu
\end{bmatrix}
$$
and $K_m$ are constant $2 \times 2$ matrices. Diagonalization of
$K_0$ is done by introducing the variable $\vec{y}=(y_1, y_2)^{\top}$ by
\begin{equation}
\begin{bmatrix}
w \\
\\
\hat{\eta}
\end{bmatrix}
=
\begin{bmatrix}
1 & \displaystyle -\frac{1}{\mathsf{g}C_{\rho}} \\
& \\
0 & 1
\end{bmatrix}
\begin{bmatrix}
y_1 \\
\\
y_2
\end{bmatrix}, \label{FS6}
\end{equation}
which reduces the equation to
\begin{equation}
s\frac{d}{ds}
\begin{bmatrix}
y_1 \\
\\
y_2
\end{bmatrix}
=\Big(
\begin{bmatrix}
0 & 0 \\
0 & -\nu
\end{bmatrix}+
\sum_{m\geq 1}s^mQ_m\Big)
\begin{bmatrix}
y_1 \\
\\
y_2
\end{bmatrix}
. \label{FS7}
\end{equation}
Applying the recipe prescribed in the proof of \cite[p.120, Chapter 4, Lemma]{CoddingtonL}, we have a fundamental matrix of solutions of \eqref{FS6}:
\begin{equation}
\Phi_{\vec{y}}=(I+\sum_s^mP_m)
\begin{bmatrix}
1 & 0 \\
0 & s^{-\nu}
\end{bmatrix}
. \label{FS8}
\end{equation}
The corresponding fundamental matrix of solutions of \eqref{FS5} is
\begin{align}
\Phi_{w,\hat{\eta}}&=
\begin{bmatrix}
1 & -\frac{1}{\mathsf{g}C_{\rho}} \\
0 & 1 
\end{bmatrix}\Phi_{\vec{y}} \nonumber \\
&=
\begin{bmatrix}
1 & \displaystyle -\frac{1}{\mathsf{g}C_{\rho}}s^{-\nu} \\
& \\
0 & s^{-\nu}
\end{bmatrix}
(I+\sum_{m\geq 1}s^mP_m), \label{FS9}
\end{align}
and the corresponding fundamental matrix of solutions of
\eqref{8a}\eqref{8b} is
\begin{align}
\mbox{\boldmath$\Phi$}_S&=
\begin{bmatrix}
1 & 0 \\
0 & s^{\nu}
\end{bmatrix}
\Phi_{w,\hat{\eta}} \nonumber \\
&=
\begin{bmatrix}
1 & \displaystyle -\frac{1}{\mathsf{g}C_{\rho}}s^{-\nu} \\
& \\
0 & 1
\end{bmatrix}
(I+\sum_{m\geq 1}s^mP_m). \label{FS10}
\end{align}
Put
\begin{equation}
\mbox{\boldmath$\Phi$}_{S}=[\mbox{\boldmath$\varphi$}_{S1} \quad \mbox{\boldmath$\varphi$}_{S2}].
\end{equation}
We have to show \eqref{FS1a}, \eqref{FS1b}. If we write
$$\sum_{m\geq 1}s^mP_m=(p_{ij})_{i,j=1,2}=(p_{ij}(s))_{i,j=1,2}, $$
we have $p_{ij}=[s]_1$, and
\begin{align*}
\mbox{\boldmath$\varphi$}_{S1}&=
\begin{bmatrix}
\displaystyle 1+p_{11}-\frac{1}{\mathsf{g}C_{\rho}}s^{-\nu}p_{21} \\
\\
p_{21}
\end{bmatrix}, \\
\mbox{\boldmath$\varphi$}_{S2}&=
\begin{bmatrix}
\displaystyle -\frac{1}{\mathsf{g}C_{\rho}}s^{-\nu}(1+p_{22})+p_{12} \\
\\
1+p_{22}
\end{bmatrix}.
\end{align*}
Therefore \eqref{FS1b} is seen. Buy \eqref{FS1a} must be proven. 

Suppose $\nu \not\in \mathbb{N}$. Suppose $p_{21}\not=0$. Then 
$p_{21}=Cs^{\mu}(1+[s]_1), C\not=0, 1\leq \mu
\in \mathbb{N}$. Then $\mu-\nu \not\in \mathbb{N}$. Since
\begin{align*}
w&=1+[s]_1-\frac{C}{\mathsf{g}C_{\rho}}s^{\mu-\nu}(1+[s]_1), \\
\hat{\eta}&=Cs^{\mu-\nu}(1+[s]_1)
\end{align*}
satisfy the equation \eqref{FS3a}, 
\begin{align*}
s\frac{dw}{ds}&=[s]_1-\frac{C(\mu-\nu)}{\mathsf{g}C_{\rho}}s^{\mu-\nu}(1+[s]_1), \\
\hat{A}_{11}w+\hat{A}_{12}\hat{\eta}&=[s]_1+
\frac{C\nu}{\mathsf{g}C_{\rho}}s^{\mu-\nu}(1+[s]_1)
\end{align*}
implies $-(\mu-\nu)=\nu$, say, $\mu=0$, a contradiction. Therefore we can calim
$p_{21}=0$ if $\nu \not\in \mathbb{N}$. 

Suppose $\nu \in \mathbb{N}$. If $p_{21}\not=0$, then 
$p_{21}=Cs^{\mu}(1+[s]_1), C\not=0, 1\leq \mu \in \mathbb{N}$. Suppose $\mu-\nu <0$. Then
\begin{align*}
s\frac{dw}{ds}& \sim -\frac{C(\mu-\nu)}{\mathsf{g}C_{\rho}}s^{\mu-\nu}, \\
\hat{A}_{11}w+\hat{A}_{12}\hat{\eta}& \sim \frac{C\nu}{\mathsf{g}C_{\rho}}s^{\mu-\nu},
\end{align*}
would imply $\mu=0$, a contradiction. Therefore $\mu -\nu \geq 0$. Suppose $\mu=\nu$.
Then 
$$s\frac{dw}{ds}=[s]_1
$$ and
$$\hat{A}_{11}w+\hat{A}_{12}\hat{\eta}=[s]_1+\frac{C\nu}{\mathsf{g}C_{\rho}},$$
a contradiction. Therefore $\mu >\nu$ or $\mu\geq \nu+1$, that is,
$p_{21}=s^{\nu+1}[s]_0$, and
\begin{align*}
&w=1+[s]_1+s[s]_0= 1+[s]_1, \\
&\eta=s^{\nu+1}[s]_0.
\end{align*}
This completes the proof. $\square$\\

Here let us note a consequence of this Proposition. Namely, since
$(z_+-z)^{-\nu}$ does not belong to $L^2(\rho dz)$ for $\nu >1$, we see that $\vec{w}$ corresponding to $\mbox{\boldmath$\varphi$}_{S2}$ does not belong to $\mathfrak{X}_l$. Therefore,
if $\lambda$ is a positive eigenvalue and $\vec{w}=(u,w)^{\top}$ is an eigenfunction associated with $\lambda$, then the corresponding $(w,\eta)^{\top}$ should coincide with 
$\alpha\mbox{\boldmath$\varphi$}_{S1}$ with a constant $\alpha$ which is $ \not=0$. Thus
$$w=\alpha+[z_+-z]_1, \qquad \alpha\not=0.$$

Note that this does not mean that $w=\alpha\mbox{\boldmath$\varphi$}_{S1}$ always ensures the boundary condition $w(0)=0$. When it is the case, $\lambda$ is an eigenvalue, that is, it is a non generic case. \\

So, we can claim
\begin{Proposition}
Let $l >0$. Any positive eigenvalue of $\vec{L}_l$ is simple.
\end{Proposition}

Anyway,
taking small positive parameter $\varepsilon$, consider the perturbation
\begin{equation}
\mbox{\boldmath$\xi$}=\varepsilon u(z)\sin(\sqrt{\lambda}t)\sin(lx)\mathbf{e}_1+
\varepsilon w(z)\sin(\sqrt{\lambda}t)\cos( lx)\mathbf{e}_3, 
\end{equation}
which is a solution of the wave equation
$$\frac{\partial^2\mbox{\boldmath$\xi$}}{\partial t^2}+\mathbf{L}\mbox{\boldmath$\xi$} =0. $$
Here $\vec{w}=(u,w)^{\top}$ is an eigenfunction of 
$$\vec{L}_l\vec{w}=\lambda\vec{w} $$
with $l > 0, lx_+\in 2\pi\mathbb{Z}$, associated with an eigenvalue $\lambda >0$. Let us normalize $\vec{w}$ so that $w = 1+[z_+-z]_1$.
This takes the value
$$\mbox{\boldmath$\xi$}=\varepsilon u(z_+)\sin(\sqrt{\lambda}t)\sin(lx)\mathbf{e}_1+
\varepsilon  \sin(\sqrt{\lambda}t)\cos( lx)\mathbf{e}_3 $$
on the vacuum boundary $\{ z=z_+\}$. Actually the trace  $u(z_+)$ exists, since 
we have
$$ u=-\frac{1}{l}\Big(\frac{\eta}{c^2\rho}+\frac{dw}{dz}\Big)=[z_+-z]_0 $$
thanks to \eqref{FS1a}. If we denote by $\underline{x}, \underline{y}, \underline{z}$ the Eulerian coordinates, then the vacuum boundary is represented as
\begin{subequations}
\begin{align}
&\underline{x}=x+\varepsilon u(z_+)\sin(\sqrt{\lambda}t)\sin(lx), \label{v1a} \\
&\underline{z}=z_++\varepsilon\sin(\sqrt{\lambda}t)\cos(lx). \label{v1b}
\end{align}
\end{subequations}
\eqref{v1a} can be solved as 
$$x=\Phi(t,\underline{x})=\underline{x}+O(\varepsilon),$$
provided that $\varepsilon$ is sufficiently small. Hence, in view of \eqref{v1b}, 
the vacuum boundary is a vibrating surface described by
\begin{align*}
\underline{z}&=z_++\varepsilon \sin(\sqrt{\lambda}t)\cos(l\Phi(t,\underline{x})) \\
&=z_++\varepsilon\sin(\sqrt{\lambda}t)\cos(l\underline{x}) +O(\varepsilon^2).
\end{align*}
Since $\varepsilon\not=0$, we observe that this surface is vibrating around $\underline{z}=z_+$.\\

This gives the configuration of standing waves. The configuration of progressive waves can be given by 
\begin{equation}
\mbox{\boldmath$\xi$}=\varepsilon u(z)\sin(lx-\sqrt{\lambda}t)\mathbf{e}_1+\varepsilon w(z)
\cos(lx-\sqrt{\lambda}t)\mathbf{e}_3.
\end{equation}
In this case $\mbox{\boldmath$\xi$}|_{t=0}\not=0$, and the initial perturbation of the Eulerian coordinate
$$\Big(\underline{\mathbf{x}}-\mathbf{x}\Big)\Big|_{t=0}=\varphi(0,\mathbf{x})-\mathbf{x}=
\mbox{\boldmath$\xi$}(0,\mathbf{x})
=\varepsilon u(z)\sin(lx)\mathbf{e}_1+\varepsilon w(z)\cos(lx)\mathbf{e}_3$$
does not identically vanish. However the wave equation \eqref{WE} and the reduction to the eigenvalue problem \eqref{4002}
are the same to the case with $\mbox{\boldmath$\xi$}|_{t=0}=0$. Actually we can assume
that $\mathrm{det}J(t, \mathbf{x}) \not=0$, since
$J(t,\mathbf{x})=I+O(\varepsilon)$, provided that $\varepsilon$ is small. 
The vacuum
boundary is represented as
\begin{subequations}
\begin{align}
\underline{x}&=x+\varepsilon u(z_+)\sin(lx-\sqrt{\lambda}t), \label{v2a} \\
\underline{z}&=z_++\varepsilon \cos(lx-\sqrt{\lambda}t). \label{v2b}
\end{align}
\end{subequations}
\eqref{v2a} can be solved as 
$$x=\Phi^{PW}(t,\underline{x})=\underline{x}+O(\varepsilon), $$
provided that $\varepsilon$ is sufficiently small, and, in view of \eqref{v2b},  the vacuum boundary is a moving surface described by
\begin{align*}
\underline{z}&=z_++\varepsilon\cos(l\Phi^{PW}(t,\underline{x})-\sqrt{\lambda}t) \\
&=z_++\varepsilon\cos(l\underline{x}-\sqrt{\lambda}t)+
O(\varepsilon^2).
\end{align*}

\section{Existence of `g-modes' and `p-modes'}

\subsection{g-modes}

We introduce 

\begin{Assumption}\label{PN}
It holds
\begin{equation}
\mathscr{N}^2(z) >0 \quad\mbox{for}\quad 0\leq z \leq z_+.
\end{equation}
\end{Assumption}

Under this assumption we are going to prove the following

\begin{Theorem}\label{Th.Exg}
Let $l > 0$ and $ lx_+ \in 2\pi\mathbb{Z}$. Suppose Assumption \ref{PN}. There exists a sequence of positive eigenvalues $(\lambda_{-n}^{(l)})_{n\in \mathbb{N}}$ of $\mathbf{L}$ 
such that $\lambda_{-n}^{(l)} \rightarrow 0$ as $ n \rightarrow \infty$. They are simple eigenvalues of $\vec{L}_l$.
\end{Theorem}

Proof. We rewrite the equation \eqref{11} as
\begin{equation}
-\frac{d}{dz}\Big(\frac{1}{\rho}\frac{d\eta}{dz}\Big)=
\Big[\frac{\lambda}{c^2}-l^2\Big(1-\frac{\mathscr{N}^2}{\lambda}\Big)\Big]
\frac{1}{\rho}\eta. \label{4.1}
\end{equation}
We consider the eigenvalue problem for the eigenvalue $\Lambda$:
\begin{equation}
-\frac{d}{dz}\Big(\frac{1}{\rho}\frac{d\eta}{dz}\Big)+
\Big(l^2-\frac{\lambda}{c^2}\Big)\frac{1}{\rho}\eta=
\Lambda \frac{l^2\mathscr{N}^2}{\rho}\eta, \label{4.2}
\end{equation}
in which $\lambda$ in the left-hand side is considered as a  parameter. If 
$\displaystyle \frac{1}{\lambda}=\Lambda$, then \eqref{4.2} coincides with \eqref{4.1}.

Taking $\lambda_0 <l\mathsf{g}$, we consider $0\leq \lambda \leq \lambda_0$. 
Actually, when a solution $\eta \in L^2(dz/c^2\rho)\subset L^2(dz/\rho)$ is given, the associated solution $w$ by which $(w,\eta)$ gives an eigenvector function should be determined by
$$
A_{21}w=\frac{l^2\mathsf{g}^2}{\lambda}\Big(1-\frac{\lambda^2}{l^2\mathsf{g}^2}\Big)\rho w = -\Big(\frac{d\eta}{dz}+A_{22}\eta\Big) =
-\Big(\frac{d\eta}{dz}+\frac{l^2\mathsf{g}}{\lambda}\eta\Big).
$$
In this point of view, we consider $\lambda$ restricted to $0\leq \lambda <l\mathsf{g}$,
that is, $A_{21}>0$.

Then for each fixed $\lambda \in [0,\lambda_0]$, we have a sequence of eigenvalues $\Lambda_n(\lambda)$ such that
$$\Lambda_1(\lambda) <\Lambda_2(\lambda)<\cdots <
\Lambda_n(\lambda) < \cdots \rightarrow +\infty.$$
 In fact we perform the Liouville transformation
\begin{align}
&\zeta=\int_0^z\sqrt{\frac{\kappa}{a}}dz, \qquad v=(a\kappa)^{\frac{1}{4}}\eta, \nonumber  \\
&q=\frac{b}{\kappa}+\frac{1}{4}\frac{a}{\kappa}\Big[\frac{d^2}{dz^2}\log(a\kappa)-
\frac{1}{4}\Big(\log(a\kappa)\Big)^2+ \nonumber \\
&+\Big(\frac{d}{dz}\log a\Big)\Big(\frac{d}{dz}\log (a\kappa)\Big)\Big]
\end{align}
for
\begin{equation}
a=\frac{1}{\rho}, \quad b=\Big(l^2-\frac{\lambda}{c^2}\Big)\frac{1}{\rho},
\qquad \kappa=\frac{l^2\mathscr{N}^2}{\rho},
\end{equation}
which transform \eqref{4.2}
to
\begin{equation}
-\frac{d^2v}{d\zeta^2}+qv=\Lambda v.
\end{equation}
We put
\begin{equation}
q=q(\zeta;\lambda)=-\frac{1}{l^2c^2\mathscr{N}^2}\lambda +q_0(\zeta)
\end{equation}
with
\begin{align}
q_0(\zeta)&=\frac{1}{\mathscr{N}^2}+
\frac{1}{4}\frac{1}{l^2\mathscr{N}^2}\Big[
\frac{d^2}{dz^2}\log\Big(\frac{\mathscr{N}^2}{\rho^2}\Big)-\frac{1}{4}
\log\Big(\frac{\mathscr{N}^2}{\rho^2}\Big)+ \nonumber \\
&+\Big(\frac{d}{dz}\log\frac{1}{\rho}\Big)\Big(\frac{d}{dz}\log
\Big(\frac{\mathscr{N}^2}{\rho^2}\Big)\Big].
\end{align}

We see that 
$$\zeta_+:=\int_0^{z_+}\sqrt{\frac{\kappa}{a}}dz=l\int_0^{z_+}\mathscr{N}dz <\infty
$$
and
$$\zeta_+-\zeta=l\mathscr{N}(z_+)(z_+-z)(1+[z_+-z]_1).
$$
Thus the interval $[0, z_+ ]$ is mapped onto $[0, \zeta_+ ]$．We see
$$q=\frac{\nu(\nu+2)}{4}\frac{1}{(\zeta_+-\zeta)^2}(1+[\zeta_+-\zeta]_1).
$$
Therefore for each fixed $\lambda \in [0, \lambda_0]$ we have simple eigenvalues
$$\Lambda_1(\lambda)<\Lambda_2(\lambda) < \cdots 
<\Lambda_n(\lambda) < \cdots \rightarrow +\infty.
$$
The $n$-th eigenvalue is given by the Max-Min principle as follows:

For any set $Y=\{y_1, \cdots, y_n\} \subset \mathfrak{E} $ we put
$$\mathfrak{d}(Y; \lambda)=\inf \{ Q[v;\lambda]\   |\   (v|y)_{\mathfrak{E}}=0 \quad\mbox{for}\quad y \in Y, v \in \mathfrak{E}_1, \|v\|_{\mathfrak{E}}=1 \}.
$$
Here $\mathfrak{E}=L^2(0,\zeta_+)$,
$$Q[v;\lambda]=\int_0^{\zeta_+}
\Big(\Big|\frac{dv}{d\zeta}\Big|^2+q(\zeta;\lambda)+K_0\Big)|v|^2d\zeta,
$$
for $$
q(\zeta; \lambda) +K_0 \geq  1+\frac{K_1}{(\zeta_+-\zeta)^2}, \qquad \frac{3}{4} < K_1 
<\frac{\nu(\nu+2)}{4}.
$$
We take $K_0$ independent of $\lambda \in [0,\lambda_0]$. We put
$$\mathfrak{E}_1=\{ v | 
\int_0^{\zeta_+}
\Big(\Big|\frac{dv}{d\zeta}\Big|^2+q_0(\zeta)+K_0\Big)|v|^2d\zeta <\infty \}.
$$
Here note that $q(\zeta;\lambda)\leq q_0(\zeta)$, for $\lambda \geq 0$, so that
$Q[v;\lambda] \leq Q[v;0]$.  The space $\mathfrak{E}_1$ is independent of $\lambda$. 
Then it is known that
$$\Lambda_n(\lambda)=\sup\{ \mathfrak{d}(Y;\lambda)\   |\   Y\subset \mathfrak{E},
\sharp Y =n\} -K_0.
$$

As for the theory of the Max-Min principle see e.g., \cite[Chapter 11]{Helffer}. The above characterization of $\Lambda_n$ is given as 
\cite[p.144, (11.3.1)]{Helffer}\\

We claim

\begin{Proposition}\label{Cg}
The function $\lambda \mapsto \Lambda_n(\lambda)$ is continuous on
$[0,\lambda_0]$.
\end{Proposition}

Proof. Let $0 \leq \lambda +\delta\lambda, \lambda \leq \lambda_0$. 
Then
$$Q[v;\lambda+\delta\lambda]-Q[v;\lambda]=
\Big(-\int_0^{\zeta_+}\frac{1}{l^2c^2\mathscr{N}^2}|v|^2d\zeta\Big)\cdot\delta\lambda.
$$
But
$$\int_0^{\zeta_+}\frac{1}{l^2c^2\mathscr{N}^2}|v|^2d\zeta \leq M
\int_0^{\zeta_+}\frac{1}{\zeta_+-\zeta}|v|^2d\zeta,
$$
since $\displaystyle c^2\sim \frac{\mathsf{g}}{\nu}(z_+-z)$ and $c^2 >0$. On the other hand
\begin{align*}
Q[v;\lambda]&=\int_0^{\zeta_+}\Big(
\Big|\frac{dv}{d\zeta}\Big|^2+\Big(-\frac{1}{l^2c^2\mathscr{N}^2}\lambda +q_0(\zeta)+K_0\Big)|v|^2d\zeta \\
&\geq \frac{1}{M'}\int_0^{\zeta_+}\frac{1}{(\zeta_+-\zeta)^2}|v|^2d\zeta,
\end{align*}
where $M'$ can be independent of $\lambda$. Therefore
$$\int_0^{\zeta_+}\frac{1}{l^2c^2\mathscr{N}^2}|v|^2d\zeta \leq \zeta_+M'Q[v;\lambda],
$$
and
$$|Q[v;\lambda+\delta\lambda]-Q[v;\lambda]|\leq
\zeta_+MM'Q[v;\lambda]|\delta\lambda|.$$
Thus we have
\begin{align*}
&-\zeta_+MM'|\delta\lambda|\mathfrak{d}(Y;\lambda_0) \leq \\
&\leq -\zeta_+MM'|\delta\lambda|\mathfrak{d}(Y; \lambda+\delta\lambda)
\leq \mathfrak{d}(Y; \lambda+\delta\lambda)-\mathfrak{d}(Y;\lambda) \\
&\leq \zeta_+MM'|\delta\lambda|\mathfrak{d}(Y; \lambda) \leq
\zeta_+MM'|\delta\lambda|\mathfrak{d}(Y;0).
\end{align*}

This estimate implies the Lipschitz continuity of $\Lambda_n(\lambda)=
\sup\mathfrak{d}(Y;\lambda)$. $\square$\\

Now, since $\Lambda_n(\lambda_0) \rightarrow +\infty$, we can find $n_0$ such that
 $\displaystyle \Lambda_n(\lambda_0) >\frac{1}{\lambda_0}$ for $n \geq n_0$. Note that
$\Lambda_n(\lambda)\geq \Lambda_n(\lambda_0)$ for $\lambda \in [0, \lambda_0]$. By Proposition \ref{Cg} the function $$f : \lambda \mapsto
\lambda -\frac{1}{\Lambda_n(\lambda)}$$
is continuous on $[0, \lambda_0]$ and
$$f(0)=-\frac{1}{\Lambda(0)} <0, \qquad f(\lambda_0)=\lambda_0-\frac{1}{\Lambda_n(\lambda_0)}>0.
$$
 Therefore there exists at least one $\lambda \in ]0, \lambda_0[$ such that
$f(\lambda)=0$, that is,
$\displaystyle \Lambda_n(\lambda)=\frac{1}{\lambda}$. Although we cannot claim that the solution is unique, we denote by $\lambda_{-n}$ one of the solutions. Then
$$\lambda_{-n}=\frac{1}{\Lambda_n(\lambda_{-n})} \leq \frac{1}{\Lambda_n(\lambda_0)}
\rightarrow 0.$$ Thus completes the proof of Theorem \ref{Th.Exg}\\

\subsection{p-modes}

In the same way we can prove the following

\begin{Theorem}
Let $l >0$ and $ lx_+ \in 2\pi\mathbb{Z}$. There exists a sequence of positive eigenvalues $(\lambda_n^{(l)})_{n \in \mathbb{N}}$ 
of $\mathbf{L}$ such that
$\lambda_n^{(l)} \rightarrow +\infty$ as $n \rightarrow \infty$. They are simple eigenvalues of $\vec{L}_l$.
\end{Theorem}

Proof. In order to deal with the equation \eqref{4.1} we consider the eigenvalue problem for the eigenvalue $\Lambda$ :
\begin{equation}
-\frac{d}{dz}\Big(\frac{1}{\rho}\frac{d\eta}{dz}\Big)+
l^2(1-\mu \mathscr{N}^2)\frac{1}{\rho}\eta=\Lambda\frac{1}{c^2\rho}\eta, \label{p1}
\end{equation}
in which $\mu$ is a parameter, which stands for $1/\lambda$. If $\mu=1/\lambda=1/\Lambda$, then \eqref{p1} coincides with \eqref{4.1} and this $\lambda$ gives a solution. 

Taking $\mu_0=1/\lambda_0 <1/l\mathsf{g}$, we consider $0\leq \mu \leq \mu_0$. We are providing $A_{21}<0$ by the restriction that $\lambda >l\mathsf{g}$. 

Then for each fixed $\mu \in [0, \mu_0]$, we have a sequence of eigenvalues $(\Lambda_n(\mu))_{n\in\mathbb{N}}$ such that
$$\Lambda_1(\mu) < \Lambda_2(\mu) < \cdots <
\Lambda_n(\mu)< \cdots \rightarrow +\infty.$$
In fact, we perform the Liouville transformation
\begin{align}
&\zeta=\int_0^z\sqrt{\frac{\kappa}{a}}dz, \qquad v=(a\kappa)^{\frac{1}{4}}\eta, \nonumber  \\
&q=\frac{b}{\kappa}+\frac{1}{4}\frac{a}{\kappa}\Big[\frac{d^2}{dz^2}\log(a\kappa)-
\frac{1}{4}\Big(\log(a\kappa)\Big)^2+ \nonumber \\
&+\Big(\frac{d}{dz}\log a\Big)\Big(\frac{d}{dz}\log (a\kappa)\Big)\Big]
\end{align}
for
\begin{equation}
a=\frac{1}{\rho}, \quad b=l^2(1-\mathscr{N}^2\mu)\frac{1}{\rho},
\qquad \kappa=\frac{1}{c^2\rho},
\end{equation}
which transform \eqref{p1}
to
\begin{equation}
-\frac{d^2v}{d\zeta^2}+qv=\Lambda v.
\end{equation}
We put 
\begin{align}
&q=q(\zeta;\mu)=-l^2c^2\mathscr{N}^2\mu +q_0(\zeta) \\
& q_0(\zeta)=l^2c^2+
\frac{1}{4}c^2\Big[
\frac{d^2}{dz^2}\log\Big(\frac{\mathscr{N}^2}{\rho^2}\Big)-\frac{1}{4}
\log\Big(\frac{\mathscr{N}^2}{\rho^2}\Big)+ \nonumber \\
&+\Big(\frac{d}{dz}\log\frac{1}{\rho}\Big)\Big(\frac{d}{dz}\log
\Big(\frac{\mathscr{N}^2}{\rho^2}\Big)\Big].
\end{align}

We see that 
$$\zeta_+:=\int_0^{z_+}\sqrt{\frac{\kappa}{a}}dz=\int_0^{z_+}\frac{1}{c}dz <\infty
$$
and
$$\zeta_+-\zeta=2\sqrt{\frac{\nu}{\mathsf{g}}}\sqrt{z_+-z}(1+[z_+-z]_1),
$$
since $$c^2=\frac{\mathsf{g}}{\nu}(z_+-z)(1+[z_+-z]_1).$$
Thus the interval $[0, z_+]$ is mapped onto $[0,\zeta_+]$．We see
$$q=\frac{(2\nu+1)(2\nu+3)}{4}\frac{1}{(\zeta_+-\zeta)^2}(1+[(\zeta_+-\zeta)^2]_1).
$$

Therefore for each fixed $\mu \in [0, \mu_0]$ we have simple eigenvalues
$$\Lambda_1(\mu)<\Lambda_2(\mu) < \cdots 
<\Lambda_n(\mu) < \cdots \rightarrow +\infty.
$$
The $n$-th eigenvalue is given by the Max-Min principle as follows:

For any set $Y=\{y_1, \cdots, y_n\} \subset \mathfrak{E} $ we put
$$\mathfrak{d}(Y; \mu)=\inf \{ Q[v;\mu]\   |\   (v|y)_{\mathfrak{E}}=0 \quad\mbox{for}\quad y \in Y, v \in \mathfrak{E}_1, \|v\|_{\mathfrak{E}}=1 \}.
$$
Here $\mathfrak{E}=L^2(0,\zeta_+)$,
$$Q[v;\mu]=\int_0^{\zeta_+}
\Big(\Big|\frac{dv}{d\zeta}\Big|^2+q(\zeta;\mu)+K_0\Big)|v|^2d\zeta,
$$
for $$
q(\zeta; \mu) + K_0\geq 1+\frac{K_1}{(\zeta_+-\zeta)^2}, \qquad \frac{3}{4} <  K_1 < \frac{2\nu+1)(2\nu+3)}{4}.
$$
We take $K_0$ independent of $\mu \in [0,\mu_0]$. We put
$$\mathfrak{E}_1=\{ v | 
\int_0^{\zeta_+}
\Big(\Big|\frac{dv}{d\zeta}\Big|^2+q_0(\zeta)+K_0\Big)|v|^2d\zeta <\infty \}.
$$
Here note that, since $|l^2c^2\mathscr{N}^2|\leq M$ uniformly, $M$ being a sufficiently large constant, we have
$$|q(\zeta;\lambda)-q_0(\zeta)|\leq M\mu_0 $$
so that
$$|Q[v;\mu]-Q[v;0]|\leq M\mu_0\|v\|_{\mathfrak{E}}^2
$$
for $0\leq \mu\leq \mu_0$.
The space $\mathfrak{E}_1$ is independent of $\mu$. 
Then it is known that
$$\Lambda_n(\mu)=\sup\{ \mathfrak{d}(Y;\mu) | Y\subset \mathfrak{E},
\sharp Y =n\} -K_0.
$$

We claim

\begin{Proposition}\label{Cp}
The function $\mu \mapsto \Lambda_n(\mu)$ is continuous on
$[0,\mu_0]$.
\end{Proposition}

Proof. Let $0\leq \mu, \mu+\delta\mu \leq \mu_0$. Then
$$Q[v;\mu+\delta\mu]-Q[v;\mu]=\Big(-\int_0^{\zeta_+}l^2c^2\mathscr{N}^2|v|^2d\zeta\Big)\cdot\delta\mu.
$$
But
$$\int_0^{\zeta_+}l^2c^2\mathscr{N}^2|v|^2d\zeta \leq M\int_0^{\zeta_+}|v|^2d\zeta =M$$
for $\|v\|_{\mathfrak{E}}=1$. Hence
$$|Q[v;\mu+\delta\mu]-Q[v;\mu]|\leq M|\delta\mu|.
$$
Easily this implies 
$$ |\Lambda_n(\mu+\delta\mu)-\Lambda_n(\mu)| \leq M|\delta\mu|.
$$
$\square$\\

Now, since $\Lambda_n(\mu_0) \rightarrow +\infty$ as $n \rightarrow \infty$, we can find $n_0$ such that $\displaystyle \Lambda_n(\mu_0) >\frac{1}{\mu_0}$ for $n \geq n_0$. Note that we can suppose that 
$\displaystyle \Lambda_n(\mu) \geq \Lambda_n(\mu_0)-M\mu_0 >\frac{1}{\mu_0}-M\mu_0 >0$ 
for $0\leq \mu \leq \mu_0$, provided that $\mu_0$ is sufficiently small.
Then the function
$$f: \mu \mapsto \mu -\frac{1}{\Lambda_n(\mu)} $$
is continuous on $[0,\mu_0]$ thanks to Proposition \ref{Cp} and 
$$f(0) =-\frac{1}{\Lambda_n(0)}<0, \qquad  f(\mu_0)=\mu_0-\frac{1}{\Lambda_n(\mu_0)} >0.$$
Therefore there exists at least one $\mu$, say
$\mu_n$, in $]0,\mu_0[$ such that $f(\mu_n)=0$, that is,
$\displaystyle \mu_n=\frac{1}{\Lambda_n(\mu_n)},$ or 
$\displaystyle \lambda_n:=\frac{1}{\mu_n}=\Lambda_n\Big(\frac{1}{\lambda_n}\Big)$,
say, $\lambda=\lambda_n$ satisfies \eqref{4.1}, and
$$\lambda_n=\Lambda_n\Big(\frac{1}{\lambda_n}\Big)\geq \Lambda_n\Big(
\frac{1}{\lambda_0}\Big)-M\mu_0\quad \rightarrow +\infty.$$
This completes the proof of the Theorem.  $\square$.\\

Note that we need not Assumption \ref{PN} in order to prove the existence of p-modes.

\section{Absence of continuous spectrum}

Let us suppose $l>0$ and consider the operator $\vec{L}_l$:
\begin{equation}
\vec{L}_l\vec{w}=
\begin{bmatrix}
L_l^u \\
\\
L_l^w
\end{bmatrix}
=\begin{bmatrix}
\displaystyle -\frac{l}{\rho}(\delta\check{P})_l \\
\\
\displaystyle \frac{1}{\rho}\frac{d}{dz}(\delta\check{P})_l+\frac{\mathsf{g}}{\rho}(\delta\check{\rho})_l
\end{bmatrix}.
\end{equation}
for $\vec{w}=(u, w)^{\top}$. 

Suppose $\lambda\not=0$. Then the equation
\begin{equation}
\vec{L}_l\vec{w}=\lambda \vec{w} \label{602}
\end{equation}
turns out to be the system \eqref{8a}\eqref{8b}, say, 
\begin{subequations}
\begin{align}
&\frac{dw}{dz}+A_{11}w+A_{12}\eta=0, \label{603a} \\
&\frac{d\eta}{dz}+A_{21}w+A_{22}\eta=0, \label{603b}
\end{align}
\end{subequations}
where 
\begin{align}
& A_{11}=-\frac{l^2\mathsf{g}}{\lambda}, \qquad A_{12}=\Big(1-\frac{l^2c^2}{\lambda}\Big)
\frac{1}{c^2{\rho}}, \nonumber  \\
& A_{21}=\frac{l^2\mathsf{g}^2}{\lambda}\Big(1-\frac{\lambda^2}{l^2\mathsf{g}^2}\Big){\rho}, \qquad A_{22}=\frac{l^2\mathsf{g}}{\lambda}.
\end{align}\\

Since $z=0$ is a regular boundary point, we can consider the solutions
$\vec{\eta}_j=(w_j, \eta_j)^{\top}, j=1,2$, such that
$$
\vec{\eta}_1=
\begin{bmatrix}
w_1 \\
\eta_1
\end{bmatrix}=
\begin{bmatrix}
0 \\
1
\end{bmatrix}
\quad\mbox{at}\quad z=0
$$ and
$$
\vec{\eta}_2=
\begin{bmatrix}
w_2 \\
\eta_2
\end{bmatrix}=
\begin{bmatrix}
1 \\
0
\end{bmatrix}
\quad\mbox{at}\quad z=0.
$$
Let us denote
$$
\mbox{\boldmath$\varphi$}_{Oj}=\vec{\eta}_j, \qquad j=1,2.
$$
Then
$$
\mbox{\boldmath$\Phi$}_O=[ \mbox{\boldmath$\varphi$}_{O1} \quad \mbox{\boldmath$\varphi$}_{O2}]
$$
is a fundamental matrix of solutions of \eqref{603a}\eqref{603b}. Of course only
$\mbox{\boldmath$\varphi$}_{O1}$ satisfies the boundary condition
$$w=0\qquad\mbox{at}\qquad z=0.
$$\\

Let us specify another fundamental matrix 
$\mbox{\boldmath$\Phi$}_S=[\mbox{\boldmath$\varphi$}_{S1} \quad \mbox{\boldmath$\varphi$}_{S2}]$ in view of asymptotic behaviors at the singular boundary $z=z_+$.
Actually, thanks to Proposition \ref{Prop.FS}, we have a fundamental matrix of solutions 
\begin{equation}
\mbox{\boldmath$\Phi$}_S =[\mbox{\boldmath$\varphi$}_{S1} \quad \mbox{\boldmath$\varphi$}_{S2}]
\end{equation}
consisting of
$\mbox{\boldmath$\varphi$}_{S1}, \mbox{\boldmath$\varphi$}_{S2}$ of the form

\begin{align}
\mbox{\boldmath$\varphi$}_{S1}&=
\begin{bmatrix}
1+O(s) \\
\\
O(s^{\nu+1})
\end{bmatrix}
, \\
\mbox{\boldmath$\varphi$}_{S2}&=
\begin{bmatrix}
\displaystyle -\frac{1}{\mathsf{g}C_{\rho}}s^{-\nu}(1+O(s)) \\
\\
1+O(s)
\end{bmatrix}
.
\end{align}

Here
\begin{equation}
s=z_+-z, \quad \Leftrightarrow z= z_+-s.
\end{equation}

Only $\mbox{\boldmath$\varphi$}_{S1}$ is admissible, since
$\eta_{S2}=1+O(s)$ does not belong to $L^2(dz/\rho)$.\\

Let 
\begin{equation}
\mbox{\boldmath$\Phi$}_S=\Phi_O\mathcal{C}(\lambda),\qquad
\mathcal{C}(\lambda)=
\begin{bmatrix}
c_{11}(\lambda) & c_{12}(\lambda) \\
c_{21}(\lambda) & c_{22}(\lambda)
\end{bmatrix}
.
\end{equation}
Here $\mathcal{C}(\lambda)$ is a non-singular matrix which is a holomorphic function of
$\lambda \in \mathbb{C}\setminus \{0\}$. See \cite[Chapter 1, Theorem 7.3]{CoddingtonL}.

Now we can claim that $\lambda$ is an eigenvalue if and only if 
\begin{equation}
D(z,\lambda):=\mathrm{det}[\mbox{\boldmath$\varphi$}_{O1}(z,\lambda) \ \mbox{\boldmath$\varphi$}_{S1}(z,\lambda)]=0.
\end{equation}
But we see 
$$D(z,\lambda)=c_{21}(\lambda)\mathrm{det}\mbox{\boldmath$\Phi$}_O(z,\lambda).
$$
Since $\mathrm{det}\mbox{\boldmath$\Phi$}_O\not=0$, we see that $\lambda$ is an eigenvalue if and only if $c_{21}(\lambda)=0$. Since $c_{21}(\lambda)$ is a holomorphic function of
$\lambda \in \mathbb{C}\setminus \{0\}$, we can claim

\begin{Theorem}
Let $l >0$.  The eigenvalues of \eqref{602} cannot accumulate to a value $\not=0 $.
\end{Theorem}

In fact, since $\vec{L}_l$ is a self-adjoint operator in $\mathfrak{X}_l$, its resolvent set contains $\mathbb{C}\setminus\mathbb{R}=\{ \lambda \  |\  \mathfrak{Im}[\lambda] \not=0 \}$ so that 
$c_{21}$ cannot vanish identically. \\

The spectrum of $\vec{L}_l$ is a subset of $[C, +\infty[ \subset \mathbb{R}$. But does it contain continuous spectra, or, does it consist of eigenvalues?  In order to examine this question, we consider $\lambda\not=0$ which is not an eigenvalue and try to solve the equation
\begin{equation}
(\vec{L}_l-\lambda)\vec{w}=\vec{f}, \label{611}
\end{equation}
where
\begin{equation}
\vec{f}=
\begin{bmatrix}
f^u \\
\\
f^w
\end{bmatrix}
\in \mathfrak{X}_l.
\end{equation}

The equation reads
\begin{subequations}
\begin{align}
&\frac{dw}{dz}+A_{11}w+A_{12}\eta=h_1, \label{613a} \\
&\frac{d\eta}{dz}+A_{21}w+A_{22}\eta=h_2, \label{613b}
\end{align}
\end{subequations}
where
\begin{equation}
\vec{h}=
\begin{bmatrix}
h_1 \\
\\
h_2
\end{bmatrix}
=
\begin{bmatrix}
\displaystyle \frac{l}{\lambda}f^u \\
\\
\displaystyle \rho\Big(f^w-\frac{l\mathsf{g}}{\lambda}f^u\Big)
\end{bmatrix}
\label{614}
\end{equation}

The solution of \eqref{613a}\eqref{613b} should be given by
\begin{equation}
\vec{\eta}(z)=
\mbox{\boldmath$\Phi$}(z)\int_0^z
\mbox{\boldmath$\Phi$}(z')^{-1}\vec{h}(z')dz'+
\mbox{\boldmath$\Phi$}(z)\mathbf{c}, \label{615}
\end{equation}
where
\begin{equation}
\mbox{\boldmath$\Phi$} = [\mbox{\boldmath$\varphi$}_{O1} \quad \mbox{\boldmath$\varphi$}_{S1}], \label{616}
\end{equation}
which is invertible since $\lambda$ is not an eigenvalue, and 
$\mathbf{c}=( c_1, c_2)^{\top}$ is a constant vector which should be chosen so that
$\vec{w}$ corresponding to $\vec{\eta}$ given by \eqref{615} belong to $ \mathsf{D}(\vec{L}_l)$.  

Actually it is possible by taking 
\begin{subequations}
\begin{align}
&c_1=-\int_0^{z_+}[\mbox{\boldmath$\Phi$}(z')^{-1}\vec{h}(z')]_1dz', \\
&c_2=0.
\end{align}
\end{subequations}
Let us show it. \\

Hereafter we suppose $\|\vec{f}\|_{\mathfrak{X}_l}\leq 1$ so that
$$\|h_1\|_{L^2(\rho dz)} \leq M_0,\qquad
\|h_2\|_{L^2(dz/\rho)}\leq M_0.$$

Let
\begin{equation}
\mathcal{C}^{-1}
=
\begin{bmatrix}
\gamma_{11} & \gamma_{12} \\
\gamma_{21} & \gamma_{22}
\end{bmatrix}.
\end{equation}

Since $$\mathrm{det}\mbox{\boldmath$\Phi$}=c_{21}\mathrm{det}\mbox{\boldmath$\Phi$}_O, $$
we have $c_{21}\not=0$ for $\lambda$ is not an eigenvalue. Therefore
$\gamma_{21}=-c_{21}/\mathrm{det}\mathcal{C} \not=0$, too. This implies that
$$\mbox{\boldmath$\varphi$}_{O1}=\gamma_{11}\mbox{\boldmath$\varphi$}_{S1}+
\gamma_{21}\mbox{\boldmath$\varphi$}_{S2} =
\gamma_{21}
\begin{bmatrix}
\displaystyle -\frac{1}{\mathsf{g}C_{\rho}}s^{-\nu}(1+O(s)) \\
\\
1+O(s)
\end{bmatrix}.
$$
Therefore we have 
\begin{equation}
\mbox{\boldmath$\Phi$}=
\begin{bmatrix}
\displaystyle -\frac{\gamma_{21}}{\mathsf{g}C_{\rho}}s^{-\nu}(1+O(s)) & 1+O(s) \\
\\
\gamma_{21}(1+O(s)) & O( s^{\nu+1})
\end{bmatrix}
\label{AsPhi}
\end{equation}
and
\begin{equation}
\mbox{\boldmath$\Phi$}^{-1}=
\begin{bmatrix}
O(s^{{\nu}+1})(1+O(s)) & \displaystyle -\frac{1}{\gamma_{21}}(1+O(s)) \\
\\
-(1+O(s)) & \displaystyle\frac{1}{\mathsf{g}C_{\rho}}s^{{\nu}}(1+O(s))
\end{bmatrix}
\label{AsPhiInverse}
\end{equation}

Then, thanks to \eqref{AsPhiInverse}, we see
\begin{align*}
&[\mbox{\boldmath$\Phi$}(z')^{-1}\vec{h}(z')]_1dz'= \\
&=\Big([\mbox{\boldmath$\Phi$}(z')^{-1}]_{11}h_1(z')+[\mbox{\boldmath$\Phi$}(z)^{-1}]_{12}h_2(z')\Big)dz' \\
&=\Big(O((s')^{{\nu}+1})h_1(z_+-s')+O(1)h_2(z_+-s')\Big)ds'.
\end{align*}
Since
\begin{align*}
&\int_0^{s_+}s^{{\nu}+1}|h_1|ds \leq
\sqrt{\int_0^{s_+}s^{{\nu}+2}ds}
\sqrt{\int_0^{s_+}|h_1|^2s^{{\nu}}ds} \leq CM_0, \\
&\int_0^{s_+}|h_2|ds
\leq
\sqrt{\int_0^{s_+}s^{{\nu}}ds}
\sqrt{\int_0^{s_+}|h_2|^2s^{-{\nu}}ds}
\leq C M_0, 
\end{align*}
we have 
$$c_1=\int_0^{z_+}
[\mbox{\boldmath$\Phi$}(z')^{-1}\vec{h}(z')]_1dz' < \infty. $$

Thus \eqref{615} reads
\begin{align*}
& w(z)=\Phi_{11}(z)H_1(z)+\Phi_{12}(z)H_2(z), \\
&\eta(z)= \Phi_{21}(z)H_1(z) + \Phi_{22}(z)H_2(z), 
\end{align*}
with
\begin{align*}
H_1(z)&=-\int_z^{z_+}
[\mbox{\boldmath$\Phi$}(z')^{-1}\vec{h}(z')]_1dz'= \\
&=
-\int_z^{z_+}
[[\mbox{\boldmath$\Phi$}(z')^{-1}]_{11}h_1(z')+
[\mbox{\boldmath$\Phi$}(z')^{-1}]_{12}h_2(z')]dz', \\
H_2(z)&=\int_0^z
[\mbox{\boldmath$\Phi$}(z')^{-1}\vec{h}(z')]_2dz'= \\
&=
\int_0^z
[[\mbox{\boldmath$\Phi$}(z')^{-1}]_{21}h_1(z')+
[\mbox{\boldmath$\Phi$}(z')^{-1}]_{22}h_2(z')]dz'.
\end{align*}

We can claim that $w, \eta$ satisfy
$ w \in L^2(\rho dz), \eta \in L^2({dz}/{\rho})$,
the boundary conditions \eqref{wBC} is satisfied and  
the corresponding $\vec{w}=(u,w)^{\top}$ belongs to $ \mathsf{D}(\vec{L}_l)$.\\

Actually, as for the boundary condition at $z=z_+$, by using \eqref{AsPhiInverse}, we can show that 
\begin{align*}
&|H_1(z)|\leq CM_0s^{\frac{\nu+1}{2}}, \\
&|H_2(z)| \leq CM_0s^{\frac{-\nu+1}{2}},
\end{align*}
therefore, by \eqref{AsPhi}, we have
\begin{align*}
&|w(z)|\leq CM_0s^{\frac{-\nu+1}{2}} \in L^2(s^{\nu}ds), \\
&|\eta(z)|\leq CM_0s^{\frac{\nu+1}{2}} \in L^2(s^{-\nu-1}ds).
\end{align*}

Here
we recall Proposition \ref{Prop.BC},
in which  the boundary conditions \eqref{wBC} is guaranteed by dint of the above estimate
of $w(z)$, and
 we recall Proposition \ref{Prop.ueta}, noting that it holds that
$$\frac{dw}{dz}+\frac{\eta}{c^2\rho} \in L^2(\rho dz) $$
if $w \in L^2(\rho dz), \eta \in L^2(dz/c^2\rho)$ and $(w, \eta)$ satisfies \eqref{613a} with 
$h_1 \in L^2(\rho dz)$,
when
$$\frac{dw}{dz}+\frac{\eta}{c^2\rho} =\frac{l^2}{\lambda}\Big(\mathsf{g}w+\frac{\eta}{\rho}\Big)
+h_1.$$
 \\

As for the boundary condition at $z=0$, we see that $H_1(z) \rightarrow c_1$ as $ z \rightarrow 0$ and $\Phi_{11}(z)=[\mbox{\boldmath$\varphi$}_{O1}]_1(z)=0$ at $z=0$, and that 
$\Phi_{12}(z)=[\mbox{\boldmath$\varphi$}_{S1}]_1(z)=O(1)$ as $z \rightarrow 0$ and
$|H_2(z)|\leq Cz^{\frac{1}{2}}$, since 
$$\mbox{\boldmath$\Phi$}^{-1}=
\begin{bmatrix}
1 & -c_{11}/c_{21} \\
& \\
0 & 1/c_{21}
\end{bmatrix}
\mbox{\boldmath$\Phi$}_O
=O(1)
$$
as $ z \rightarrow 0$ and $(h_1, h_2) \in L^2(dz)$; Thus $w=\Phi_{11}H_1+\Phi_{12}H_2
\rightarrow 0$ as $ z \rightarrow 0$. \\

Summing up, we can claim

\begin{Theorem}
Let $l > 0$.  Let $\lambda \not=0$, and $\lambda$ is not an eigenvalue. Then for any $\vec{f}\in \mathfrak{X}_l$ the equation
$$(\vec{L}_l-\lambda)\vec{w}=\vec{f}$$ admits a solution $\vec{w} \in \mathsf{D}(\vec{L}_l)$
and
$$\|\vec{w}\|_{\mathfrak{X}_l} \leq C\|\vec{f}\|_{\mathfrak{X}_l},
$$
that is, $(\vec{L}_l-\lambda)^{-1}$ is bounded in $\mathfrak{X}_l$, or,
$\lambda$ belongs to the resolvent set of the operator $\vec{L}_l$. Therefore the spectrum of $\vec{L}_l$
consists of countable many eigenvalues, which are simple if $ \not= 0$.
\end{Theorem}

As a corollary we can claim

\begin{Theorem}
The eigenfunctions of $\vec{L}_l$ form a complete orthogonal system of $\mathfrak{X}_l$.
\end{Theorem}

Here, when $l=0$, the assertion should mean that the eigenfunctions of $L_0^w$, which are associated with the eigenvalues $(\lambda_n^{(0)})_{n=1,2,\cdots}$ (Theorem \ref{Efl0} ), is complete in $L^2([0, z_+]; \rho(z) dz)$,
by neglecting $u$ to be $0$, while we consider $\mathfrak{X}_0=\{ 0\} \times L^2([0,z_+];\rho(z)dz)$. \\

For a proof of the completeness of eigenfunctions, see \cite[p.905, X.3.4 Theorem]{DunfordS}, which can be applied 
to the unbounded self-adjoint operator $\vec{L}_l$ thanks to \cite[p.177, Chapter III, Theorem 6.15]{Kato}.

\section{Concluding remarks}

Under the Assumption \ref{PN} that $\mathscr{N}^2 \geq 1/C >0$, we have given a mathematically rigorous proof of the existence of  a sequence of eigenvalues which accumulates to $0$, say, the existence of the g-modes. However the problem to clarify the structure of the spectrum of the operator $\mathbf{L}$ for this case cannot be said to be completely solved, say, the question whether the spectrum of the self-adjoint operator $\mathbf{L}$ considered in a suitable dense subspace $\mathfrak{F}$ of $\mathfrak{H}$ can be exhausted by eigenvalues, which accumulate to $0$ and $+\infty$ or not is still open. In this sense the present result has  not yet caught up the clear conclusion for the case in which $\mathscr{N}^2=0$ everywhere given by \cite{JJTM}.

On the other hand our argument on the existence of g-modes cannot work if $\mathscr{N}^2$ takes negative values somewhere. However in the theory of astroseismology by astrophysicists the square of the Brunt-V\"{a}is\"{a}l\"{a} frequency $\mathscr{N}^2$ turns out to be negative near the surface in many realistic stellar models. (See, e.g.,  \cite{LedouxS}.) In this sense we should continue the study for the case in which the Assumption \ref{PN} does not hold.   

Of course we should not forget that this article concerns the solutions of the  linearized approximation. The ultimate aim of the mathematical study is to construct true solutions of the original non-linear equations for which the constructed time-periodic  oscillations of the linearized equations turn out to be the first approximations of the true solutions. As for barotropic and spherically symmetric evolution of the structure of the atmosphere or the self-gravitating gaseous stars, this task has been done successfully in \cite{TM.FE}, \cite{TM.OJM}, \cite{JJ.APDE}, but there remain open problems of this task for not barotropic, or, not spherically symmetric perturbations. 

Therefore even if we consider the problem of gaseous adiabatic oscillations under gravitation in the most simple situation on the flat Earth around a stratified back ground density distribution,
there remain many interesting mathematical problems still open. Mathematical difficulty
arises from the treatise of the free boundary of the gas which touches the vacuum.

\vspace{15mm}

{\bf\Large Acknowledgment}\\

This work is supported by JSPS KAKENHI Grant Number JP18K03371. The idea of the proof of the existence of g-modes was obtained during the stay of the author at the Department of Mathematics, National University of Singapore in March 4-11, 2020. The author expresses his sincerely deep  thanks to Professor Shih-Hsien Yu for the invitation and stimulating discussions, and to the Department of Mathematics, National University of Singapore for the hospitality and the financial support. The author expresses his sincerely deep thanks to the anonymous referee for the careful examination of the original manuscript and the provision of accurate comments which were inevitable for the revision of the presentation of this work.\\

{\bf\Large  The data availability statement}\\

The data that supports the findings of this study are available within the article and the 
list of  references.

\end{document}